\documentclass[11pt,reqno]{amsart}

\usepackage[numbers,sort&compress]{natbib}
\usepackage[colorlinks,citecolor=blue]{hyperref}
\usepackage{amssymb}

\usepackage{enumerate}

\usepackage{graphicx}
\usepackage{booktabs}
\usepackage{mathrsfs}
\usepackage{algorithm}
\usepackage{algorithmic}
\usepackage{graphicx}
\usepackage{subfigure}
\usepackage{epstopdf}
\usepackage{amssymb}
\usepackage{caption}
\usepackage{subfigure}
\usepackage{verbatim}
\usepackage{color}
\usepackage{amsmath}
\usepackage{cases}
\usepackage{float}
\makeatletter\@addtoreset{equation}{section}

\makeatletter
\@namedef{subjclassname@2010}{%
\textup{2010} Mathematics Subject Classification}
\makeatother

\newtheorem{thm}{Theorem}[section]

\theoremstyle{definition}

\newtheorem{rem}[thm]{Remark}

\numberwithin{equation}{section}

\makeatletter\@addtoreset{equation}{section}

\newenvironment{proofed}[1]{\par \textbf{Proof}\quad #1}{\hfill \textbf{} $\Box$ }

\newtheorem{example}{\bf{Example}}[section]

\numberwithin{equation}{section}
\numberwithin{table}{section}
\numberwithin{figure}{section}
\newcommand{\ba}{\begin{array}}\newcommand{\ea}{\end{array}}
\newcommand{\be}{\begin{eqnarray}}\newcommand{\ee}{\end{eqnarray}}
\newcommand{\beq}{\begin{equation*}}\newcommand{\eeq}{\end{equation*}}
\newcommand{\bex}{\begin{eqnarray*}}
\newcommand{\eex}{\end{eqnarray*}}
\newcommand{\bse}{\begin{subequations}}
\newcommand{\ese}{\end{subequations}}
\newcommand{\bal}{\begin{align}}
\newcommand{\eal}{\end{align}}

\usepackage{lineno}
\usepackage{enumerate}
\linenumbers
\pagewiselinenumbers
\font\tenbi=cmmib10   at 11 pt \font\sevenbi=cmmib10 at 9pt
\font\fivebi=cmmib7 at 6pt
\newfam\bifam
\textfont\bifam=\tenbi
\scriptfont\bifam=\sevenbi
\scriptscriptfont\bifam=\fivebi

\def\bi{\fam\bifam\tenbi}
\font\tendb=msbm10 at 12 pt \font\sevendb=msbm7
\newfam\dbfam
\textfont\dbfam=\tendb \scriptfont\dbfam=\sevendb

\def\r{\color{red}}

\def\u{{\bi u}}

\begin{document}
\graphicspath{{figure/};}
\baselineskip=17pt
\nolinenumbers 
\title
[New schemes for the Cahn-Hilliard-Brinkman system]
{New highly efficient and accurate numerical  scheme  for the Cahn-Hilliard-Brinkman system$^\ast$}

\author[D.W. Chen]{Dawei Chen$^{1}$}
\author[Q.Z. Ren]{Qinzhen Ren$^{1}$}
\author[M.H. Li]{Minghui Li$^{1,\dag}$}
\thanks{\hskip -12pt
$^\ast$ This research is partially supported by Scientific Research Foundation of Jiangxi Province (No.GJJ2201360) and PhD Scientific Research Foundation of Jiangxi Science and Technology Normal University (No. 2022BSQD15).\\
$^{1}$School of Mathematical Sciences, Jiangxi Science and Technology Normal University, 330038 Nanchang, China.
\\${}^{\dag}$Corresponding author.\\
Emails:  cdwstnu@163.com(D. Chen); qzhen2017@163.com(Q. Ren); minghuili\_xmu@163.com (M. Li).}
\date{}

\begin{abstract}
In this paper, based on a generalized scalar auxiliary variable approach with relaxation (R-GSAV), we construct a class of high-order backward differentiation formula (BDF) schemes with variable time steps for the Cahn-Hilliard-Brinkman(CHB) system. In theory, it is strictly proved that the designed schemes are unconditionally energy-stable. With the delicate treatment of adaptive strategies, we propose several adaptive time-step algorithms to enhance the robustness of the schemes.  More importantly, a novel hybrid-order adaptive time steps algorithm performs outstanding for the coupled system. The hybrid-order algorithm inherits the advantages of some traditional high-order  BDF adaptive strategies. A comprehensive comparison with some adaptive time-step algorithms is given, and the advantages of the new adaptive time-step algorithms
are emphasized.  Finally, the effectiveness and accuracy of the new methods are validated through a series of numerical experiments.
\end{abstract}

\maketitle

\section{Introduction}
The Cahn-Hilliard-Brinkman equation was first proposed to describe the phase separation of binary fluids in a porous medium \cite{ngamsaad2010theoretical}. This well-posed nonlinear system consists of a convective Cahn-Hilliard equation that mainly governs the interface evolution and the Brinkman equation that controls the fluid flows in complex media. The fourth-order convective Cahn-Hilliard equation features a gradient energy term that sharpens fluid interfaces and a diffusion term that smooths them. The Brinkman equation is an extension to the traditional form of Darcy's law and introduces an additional coupled term to account for diffuse interface surface tension force\cite{collins2013efficient}.
 More recently, using the convective Cahn-Hilliard equation for the evolution of the tumour and the Brinkman-Stokes type law for the fluid velocity, the CHB model has been adopted as a comprehensive mathematical model to provide deeper insights into tumour growth dynamics, understand key mechanisms and design new treatment strategies\cite{ebenbeck2019cahn,knopf2022existence,colli2023cahn,ebenbeck2021cahn,garcke2016cahn}, but computationally, it is not easy to resolve. The hybrid system is closely related to the Cahn-Hilliard-Navier-Stokes systems \cite{han2015second,chen2020novel,li2020sav,li2021new,liu2015decoupled} and many others, which are powerful mathematical models employed to simulate the behavior of multiphase flows, especially those involving immiscible fluids.
 The model also satisfies a thermodynamically consistent energy dissipation
law. Constructing schemes that preserve discrete energy dissipation laws is crucial to prevent significant numerical errors during long-time simulations\cite{Beckerrohl2008,yamaleev2009systematic,shen2010numerical,fjordholm2011well,wang2011energy,chen2013time}.
This study aims to develop an efficient numerical scheme and an adaptive time-stepping algorithm for the CHB system while discretely upholding these energy dissipation laws.

	Several works have been devoted to numerical approximation for the CHB system.  In \cite{collins2013efficient}, an energy-stable scheme was first designed for the CHB system using the convex splitting method. Combined with the locally discontinuous Galerkin method, Guo et al.
 \cite{guo2014efficient} proposed a fully discrete energy stable scheme. The convex splitting method does not easily extend to a second-order accuracy in time for the time derivative in the phase equation. Inspired by the generalized scalar auxiliary variable approach (GSAV) \cite{huang2022new}, several higher-order unconditionally energy-stable schemes are constructed in \cite{zheng2022new}. Following the ideal of \cite{zheng2022new}, Jiang et al.
	 \cite{jiang2024highly}  adopt the implicit–explicit Crank-Nicolson SAV method to design a second-order adaptive time steps scheme(IMEX-CNA)  for CHB system. The adaptive algorithm employs the first-order numerical solution for time step control, leading to high computational costs.

 To observe the long-term behavior of interface structure, it is important to design robust adaptive algorithms for reducing costs\cite{zhang2013adaptive,cheng2018multiple}. Recently, some variable-time-step  BDF energy stable schemes have been applied to deal with the gradient flow \cite{hou2023implicit,liao2022adaptive}, in which
the modified discrete energy keeps unconditional stable only when the ratio of adjacent time steps ratio $\gamma^{*}=\tau_{n+1}/\tau_{n}\leq4.8645$ with the kernel recombination technique \cite{liao2022adaptive} or some essential inequalities \cite{hou2023implicit}. However, adopting the techniques similar with \cite{hou2023implicit,liao2022adaptive} to the CHB model proves to be challenging due to the presence of a coupling term in the phase equation. Additionally, the high-order nonuniform approximation of the time derivative in the phase equation further complicates the discrete energy stability analysis. In \cite{huang2020highly}, a class of efficient, high-order, unconditionally energy-stable adaptive time-step BDF schemes were designed for gradient flow by adopting a suitable error indicator and first-order discretizations in the dynamical equation. This approach appears to offer a feasible method for developing high-order, energy-stable adaptive time-stepping algorithms for highly coupled equations. To achieve better performance, the time-adaptive algorithm in \cite{huang2020highly} partially relies on the original energy instead of the modified energy to improve the long-time accuracy of the numerical scheme. However, this approach sacrifices the unconditional energy stability in the theoretical analysis.
\par

GSAV is a modified version of the explicit-implicit BDF\(k\) schemes for general dissipative systems \cite{huang2020highly}. It enables rigorous, unified error analysis for \(k = 1,2,3,4,5\) with fixed time steps \cite{huang2022new}. Moreover, the R-GSAV method in \cite{zhang2022generalized} preserves all the advantages of the GSAV approach while directly dissipating modified energy, which is closely related to the original free energy. However, there is a scarcity of studies that introduce adaptive time-stepping strategies for the GSAV or R-GSAV methods. Further investigations indicate the adaptive strategies in \cite{huang2020highly} will not directly fit the popular  R-GSAV\cite{zhang2022generalized}. It is natural, then, to ask whether it is possible to come up with a more efficient and energy-stable adaptive time-step algorithm within the R-GSAV framework for solving highly coupled phase-field models.

This paper proposes easy-to-implement, high-order accurate, time-step adaptive, and energy-stable schemes for the CHB model. First, inspired by R-GSAV \cite{zhang2022generalized}, we introduce an extra energy dissipation equation and design a class of relaxation implicit-explicit \(k\)-order BDF (R-IMEX-BDF\(k\)) schemes with variable time steps for the CHB system. These schemes retain the advantages of the GSAV approach \cite{zheng2022new} and maintain a modified energy dissipation linked directly to the original free energy. Next, we improve the adaptive algorithms in \cite{huang2020highly} and construct the R-IMEX-BDF\(k\) schemes with the adaptive algorithm (R-IMEX-BDF\(k\)A). Finally, based on R-IMEX-BDF2A and R-IMEX-BDF3A methods, we develop a hybrid-order adaptive time-stepping algorithm for the CHB model, which significantly reduces computational costs and improves efficiency compared to the schemes with fixed time steps.\par

 The structure of the article is as follows. In Section 2, we present the R-IMEX-BDF\(k\) (\(k=1,2,3,4\)) schemes with variable time steps and prove that they keep unconditionally energy stable. Section 3 focuses on the construction and analysis of some adaptive time-stepping algorithms for the R-IMEX-BDFk schemes. Specifically, we develop the R-IMEX-BDF2A and R-IMEX-BDF3A methods, as well as a hybrid-order adaptive time-stepping algorithm tailored for the CHB model. In Section 4, we present several numerical experiments to demonstrate the accuracy and efficiency of the proposed schemes. Finally, concluding remarks are provided in Section 5.

\section{Governing equations and numerical schemes}
\setcounter{equation}{0}
\subsection{Governing equations}
We are interested in numerically solving the CHB model\cite{collins2013efficient}:\begin{subequations}\label{chb}
	\begin{align}
			&\phi_{t}=\nabla\cdot(M(\phi)\nabla\mu)-\nabla\cdot(\u\phi),
			\qquad\qquad\quad\; \label{problem1}\\
			&\mu=-\varepsilon^{2}\Delta\phi+F' (\phi),
			\qquad\qquad\qquad\;\;\; \label{problem2}\\
			&-\nabla\cdot[\nu(\phi)D(\u)]+\eta(\phi)\u=-\nabla p-\gamma\phi\nabla\mu,
			\label{problem3}\\
			&\nabla \cdot \u=0,
	\qquad\qquad\qquad\qquad\qquad\qquad\qquad\:\:\:  \label{problem4}
	\end{align}
\end{subequations}
in the domain $\Omega\subset R^{d}\ (d=2,3)$,	where $\phi$, $\mu$, $\u$, and $p$ represent the phase function, chemical potential velocity, and pressure, respectively. $\varepsilon$ is the thicknesses of transition layer, $\gamma$ is the surface tension parameter.\ $M(\phi)$, $\nu(\phi)$, and $\eta(\phi)>$0, which indicate the mobility, the fluid viscosity, and the fluid permeability, respectively.
	$F(\phi)=\frac{1}{4}(\phi^2-1)^2$ is the double-well potential.\ $D(\u):=\nabla \u+\nabla \u^T$ is the stress tensor. Here, we equip the system with periodic boundary conditions.\par
	Taking the inner products of \eqref{problem1} with $\mu$, \eqref{problem2} with $\phi_{t}$, \eqref{problem3} with $\frac{1}{\gamma}\u$ as follow:
		\begin{subequations}\label{chbE}
		\begin{align}
			&(\phi_{t},\mu)=(\nabla\cdot(M(\phi)\nabla\mu)-\nabla\cdot(\u\phi),\mu),\label{problem5}\\
			&(\mu,\phi_{t})=(-\varepsilon^{2}\Delta\phi+F' (\phi),\phi_{t}),\label{problem6}\\
			&(-\nabla\cdot[\nu(\phi)D(\u)]+\eta(\phi)\u,\frac{1}{\gamma}\u)=(-\nabla p-\gamma\phi\nabla\mu,\frac{1}{\gamma}\u)\label{problem7},
		\end{align}
	\end{subequations}
	and then adding \eqref{problem5}, \eqref{problem6} and \eqref{problem7} together, we obtain the  following energy dissipation law:
	\begin{equation}\label{dis1}
		\frac{d(E(\phi))}{dt}=-\int_{_{\Omega}}M(\phi)|\nabla\mu|^2+ \frac{1}{\gamma}\eta(\phi)|\u|^2+\frac{1}{2\gamma}\nu(\phi)|D(\u)|^2~dx,
	\end{equation}
	where $E(\phi)=\int_{_{\Omega}}(\frac{\varepsilon^{2}}{2}|\nabla\phi|^2+F(\phi))dx$ is the free energy functional.
\subsection{R-IMEX-BDFk  with variable time step}
Inspired by the  R-GSAV  approach with constant time steps  \cite{zhang2022generalized}, let $r(t)=E_1(\phi)$ and $E_1(\phi)=E(\phi)+C_0$, where constant $C_0\ge 0$, we construct the R-IMEX-BDFk scheme with variable time steps for \eqref{chb} and \eqref{dis1} below, inheriting the advantages of the IMEX-type SAV scheme with constant time steps \cite{zheng2022new} and improving its accuracy:\par

\textbf{Scheme 1 }(R-IMEX-BDF$k$):Given initial conditions $(\phi^n,\u^n,p^n,r^n)$, we update $(\phi^{n+1},\u^{n+1},p^{n+1},\\r^{n+1},\xi_k^{n+1},\zeta_k^{n+1})$\ by using the following two steps:\par
	 \textbf{Step1:}Compute an intermediate solution $(\phi^{n+1},\u^{n+1},\tilde r^{n+1})$ by using the IMEX type SAV scheme:
	 \bse\label{sche2}
	 \bal
	 &\frac{\alpha_{k}\tilde\phi^{n+1}-A_{k}({\phi^n})}{\tau_{n+1}}=\nabla\cdot(M(B_{k}(\phi^{n}))\nabla\tilde\mu^{n+1})-\nabla\cdot(B_{k}(\u^{n})B_{k}(\phi^{n}))),\label{A21}\\
	 &\tilde\mu^{n+1}=-\varepsilon^{2}\Delta\tilde\phi^{n+1}+S\tilde\phi^{n+1}+\tilde F' [(B_{k}(\phi^{n})],\label{A22}\\
	 &-\nabla\cdot[\nu( B_{k}(\phi^{n}))D(\tilde\u^{n+1})]+\eta(B_{k}(\phi^n))\tilde\u^{n+1}=-\nabla \tilde p^{n+1}-\gamma(B_{k}(\phi^n))\nabla(B_{k}(\mu^n)),\label{A23}\\
	 &\nabla\cdot\tilde\u^{n+1}=0,\label{A24} \\
	 &\frac{\tilde r^{n+1}-{r}^n}{\tau_{n+1}}=-\frac{\tilde r^{n+1}}{E_1(\tilde\phi^{n+1})}\kappa(\tilde\phi^{n+1},\tilde\u^{n+1}),\label{A25}\\
	 &\xi_k^{n+1}=\frac{\tilde r^{n+1}}{E_1(\tilde\phi^{n+1})},\label{A26}\\
&\phi^{n+1}=\zeta_k^{n+1}\tilde\phi^{n+1},\u^{n+1}=\zeta_k^{n+1}\tilde\u^{n+1},p^{n+1}=\zeta_k^{n+1}\tilde p^{n+1}\ with\ \zeta_k^{n+1}=1-(1-\xi_k^{n+1})^{k+1}\label{A27}.
	 \end{align}
	 \ese
    where $\tau_{n+1}$ is variable time step, $\tilde F' (\phi)=\phi^{3}-(S+1)\phi \ $  and $S\ge 0$ is a suitable stabilization parameter\cite{shenScalar2018a},
    $$\kappa(\phi,\u)=\int_{_{\Omega}}M(\phi)|\nabla\mu|^2+ \frac{1}{\gamma}\eta(\phi)|\u|^2+\frac{1}{2\gamma}\nu(\phi)|D(\u)|^2)~dx,$$
 The corresponding   $\ \alpha_{k}\ , A_{k}\ and\ B_{k}$ are set as  with $k$th-order extrapolation can be derived by Taylor expansion \cite{huang2020highly}, for example:\\
Second-order($k=2$):
	 \begin{eqnarray*}	 	
	 	&&\alpha_{2}=\frac{-2t^{n+1}+t^{n}+t^{n-1}}{-t^{n+1}+t^{n-1}},\\
        &&A_{2}(\phi^{n})=\frac{-t^{n+1}+t^{n-1}}{t^{n}-t^{n-1}}\phi^{n}+\frac{(-t^{n+1}+t^{n})^{2}}{(t^{n}-t^{n-1})(-t^{n+1}+t^{n-1})}\phi^{n-1},\\	 	
	 &&B_{2}(\phi^{n})=-\frac{-t^{n+1}+t^{n-1}}{t^{n}-t^{n-1}}\phi^{n}+\frac{-t^{n+1}+t^{n}}{t^{n}-t^{n-1}}\phi^{n-1}.
	 \end{eqnarray*}
	 Third-order($k=3$):\\
	 \begin{eqnarray*}
&&\alpha_{3}=(-t_{n+1}+t_{n})\left( \frac{1}{-t_{n+1}+t_{n-1}}+\frac{1}{-t_{n+1}+t_{n-2}}\right)+1,\\
&&A_{3}(\phi^{n})=\frac{(-t_{n+1}+t_{n-1})(-t_{n+1}+t_{n-2})}{(t_{n}-t_{n-1})(t_{n}-t_{n-2})}\phi^{n}-\frac{(-t_{n+1}+t_{n})(-t_{n+1}+t_{n-2})}{(-t_{n+1}+t_{n-1})(t_{n+1}-t_{n-1})(t_{n-1}-t_{n-2})}\phi^{n-1}\\
&&\qquad\qquad\ +\frac{(-t_{n+1}+t_{n})^{2}}{(-t_{n+1}+t_{n-2})(t_{n}-t_{n-2})(t_{n+1}-t_{n-2})}\phi^{n-2},\\
&&B_{3}(\phi^{n})=\frac{(-t_{n+1}+t_{n-1})(-t_{n+1}+t_{n-2})}{(t_{n}-t_{n-1})(t_{n}+t_{n-2})}\phi^{n}-\frac{(-t_{n+1}+t_{n})(-t_{n+1}+t_{n-2})}{(t_{n}-t_{n-1})(t_{n-1}+t_{n-2})}\phi^{n-1}\\
&&\qquad\qquad\ +\frac{(-t_{n+1}+t_{n})(-t_{n+1}+t_{n-1})}{(t_{n}+t_{n-2})(t_{n-1}-t_{n-2})}\phi^{n-2}.
\end{eqnarray*}	
	 \textbf{Step2:}Compute scalar auxiliary variable $r^{n+1}$ by using the following relaxation:
	 \be\label{A29}
	 r^{n+1}=\sigma^{n+1}_{0} \tilde r^{n+1}+(1-\sigma^{n+1}_{0})E_1(\phi^{n+1}), \ \sigma^{n+1}_{0}\in \Gamma,
	 \ee
	 where
	 \be\label{A30}
	 \Gamma=\left\{ \sigma\in[0,1]\ s.t.\ \frac{ r^{n+1}-\tilde{r}^{n+1}}{\tau_{n+1}}=-\delta^{n+1}\kappa(\phi^{n+1},\u^{n+1})+\frac{\tilde r^{n+1}}{E_1(\tilde\phi^{n+1})}\kappa(\tilde\phi^{n+1},\tilde{\u}^{n+1}) \right\},
	 \ee
	 with $\delta^{n+1}\geqslant0$ to be determined. The following theorem summarizes the selection of $\sigma^{n+1}_{0}\ and\ \delta^{n+1}$.

  \begin{thm}\label{chbR5}
	We choose $ \sigma^{n+1}_{0} in \ \eqref{A29} \ and \ \delta^{n+1}\  in\  \eqref{A30}$ as follows:
\begin{eqnarray*}\label{chbR6}
\begin{array}{ll}
	1.\text{If}\ \tilde r^{n+1}=E_1(\phi^{n+1}), \text{we obtain}\ \sigma^{n+1}_{0}=0\ and\\
   \centerline{$  \  \delta^{n+1}=\frac{\tilde r^{n+1}\kappa(\tilde\phi^{n+1},\tilde{\u}^{n+1})}{E_1(\tilde\phi^{n+1})\kappa(\phi^{n+1},{\u}^{n+1})}.$}
   \\
	2.\text{If}\ \tilde r^{n+1}>E_1(\phi^{n+1}),\text{we obtain}\ \sigma^{n+1}_{0}=0\ and \\
       \centerline{$ \  \delta^{n+1}=\frac{\tilde r^{n+1}-E_1(\phi^{n+1})}{\tau_{n+1}\kappa(\phi^{n+1},{\u}^{n+1})}+\frac{\tilde r^{n+1}\kappa(\tilde\phi^{n+1},\tilde{\u}^{n+1})}{E_1(\tilde\phi^{n+1})\kappa(\phi^{n+1},{\u}^{n+1})}.$}\\
	3.\text{If}\ \tilde r^{n+1}<E_1(\phi^{n+1})\ and\\
    \centerline{$\ \tilde r^{n+1}-E_1(\phi^{n+1})+\tau_{n+1}\frac{\tilde r^{n+1}}{E_1(\tilde\phi^{n+1})}\kappa(\tilde\phi^{n+1},\tilde{\u}^{n+1})\geq0,$} \\
    \ \ \text{we obtain}\ \sigma^{n+1}_{0}=0\ and  \ \delta^{n+1}\text{is the same as  case 2}.\\
	4.\text{If}\ \tilde r^{n+1}<E_1(\phi^{n+1})\ and\\
    \centerline{$\ \tilde r^{n+1}-E_1(\phi^{n+1})+\tau_{n+1}\frac{\tilde r^{n+1}}{E_1(\tilde\phi^{n+1})} \kappa(\tilde\phi^{n+1},\tilde{\u}^{n+1})<0,$}\\
 \  \ \text{we obtain}\\
 \centerline{$\ \sigma^{n+1}_{0}=1- \frac{\tau_{n+1}\tilde r^{n+1}\kappa(\tilde\phi^{n+1},\tilde{\u}^{n+1})}{E_1(\tilde\phi^{n+1})(E_1(\phi^{n+1})-\tilde r^{n+1})},$}\\
     \ \  and \ \delta^{n+1}=0.
\end{array}
\end{eqnarray*}
Then, $\sigma^{n+1}_{0}\subset\varGamma$. Given $r^{n}\geqslant 0$, we have $r^{n+1}\geqslant 0,\  \xi_k^{n+1}\geqslant 0$ in the scheme \eqref{sche2}-\eqref{A29}  and\ $r^{n+1}$ is unconditionally energy stable in the sense that
\begin{eqnarray}\label{proo1}
	r^{n+1}-{r}^n=-\tau_{n+1}\delta^{n+1}\kappa(\phi^{n+1},\u^{n+1})\leq 0, \  \ \forall n\geq0.
\end{eqnarray}
Furthermore, we obtain
\begin{eqnarray}\label{proo2}
	r^{n+1}\leq E_1(\phi^{n+1}),\  \ \forall n\geq0.
\end{eqnarray}
\end{thm}
\begin{proofed}	
	Given $\ r^{n}\geq0\ $, since $\ E_1(\tilde\phi^{n+1})\geq0\ $, it  follows from  \eqref{A25} that
	\begin{eqnarray}\label{proo4}
		\tilde r^{n+1}=\frac{r^{n}}{1+\tau_{n+1}\frac{\kappa(\tilde\phi^{n+1},\tilde\u^{n+1})}{E_1(\tilde\phi^{n+1})}}\geq0.
	\end{eqnarray}
	From \eqref{A26}, we deduce that $\xi_k^{n+1} \geq 0$, and from \eqref{A29}, we conclude that $r^{n+1} \geq 0$. By combining \eqref{A25} with \eqref{A30}, we derive \eqref{proo1}. Furthermore, combining \eqref{A29} with the above four cases, we obtain \eqref{proo2}.
	 \end{proofed}
\begin{rem}
For  $C_0>0$, the functional $E_1(\phi)$ also satisfies the energy dissipation law \eqref{dis1}. Moreover, $H^1$ stability of the numerical solution can be obtained by a similar analysis in \cite{zhang2022generalized}. Let linear positive  definite operator $L=\varepsilon^{2}\Delta+2C_0 I$, then \ $E_1(\phi)=\frac{1}{2}(L\phi,\phi)+\int_{_{\Omega}}F(\phi)dx$.
	Denote $M := r^{0} = E_1[\phi(\cdot, 0)]$. Combining \eqref{proo1} and \eqref{proo4}, we obtain $\tilde{r}^{n+1} \leq M$.  From \eqref{A27}, there
exists $C_k>0$ such that	\begin{eqnarray}\label{proo5}
		|\zeta_k^{n+1}|\leq |1-(1-\xi_k^{n+1})^{k+1}| \leq C_k \xi_k^{n+1}=C_k\frac{\tilde r^{n+1}}{E_1(\tilde\phi^{n+1})}\leq\frac{C_k M}{E_1(\tilde\phi^{n+1})}.
	\end{eqnarray}
	 With $\phi^{n+1}=\zeta_k^{n+1}\tilde\phi^{n+1}$ ,we get
	 \begin{eqnarray}\label{proo7}
	 	(L\phi^{n+1},\phi^{n+1})=(\zeta_k^{n+1})^{2}(L\tilde\phi^{n+1},\tilde\phi^{n+1})\leq\left(\frac{C_k M}{E_1(\tilde\phi^{n+1})}\right)^{2}(L\tilde\phi^{n+1},\tilde\phi^{n+1})^{2}\leq (C_k M)^2.	
	 \end{eqnarray}
\end{rem}
\section{R-IMEX-BDF\texorpdfstring{$k$}{k} method with adaptive time-stepping}
To obtain efficient numerical results, smaller time steps should be used when the energy and numerical solution of the CHB model vary dramatically, while larger time steps are suitable for smaller changes. Therefore, an adaptive time-stepping algorithm is essential for obtaining efficient numerical solutions of the CHB model. Moreover, unconditionally energy-stable schemes provide greater flexibility in designing adaptive time-stepping algorithms than conditionally energy-stable schemes, which are limited by stability constraints \cite{huang2020highly}.\par

To design an adaptive strategy, a popular approach is using an energy-based time-stepping control strategy. This involves using a functional of \( E'(t) \) to adjust the time steps for unconditionally energy-stable schemes \cite{zhang2013adaptive, cheng2018multiple, hou2023implicit, hou2023linear}.
When the energy of gradient flow models changes drastically, $E'(t)$  becomes large, prompting the use of smaller time steps. Additionally, an effective error indicator is crucial in adaptive strategy design. When the error indicator is large, smaller time steps should be employed, while larger steps are used when the indicator is small. A posterior error estimate, particularly the optimal posterior error estimate, is a good choice for this purpose \cite{chen2024recovery,feng2008posteriori}.
However, for the nonlinear CHB model, obtaining an appropriate posterior error estimate is challenging.  Fortunately, following the approach for gradient flows in \cite{huang2020highly}, to achieve accurate results, $\xi_k^{n+1}$ has to be close to 1, and $|1-\xi_k^{n+1}|$ would serve as an error indicator for the adaptive strategy. The coupled CHB system is more complex than the Cahn-Hilliard or Allen-Cahn models in \cite{huang2020highly}, and we incorporate a relaxation term in the R-IMEX-BDF$k$ scheme to improve accuracy. Therefore, $|1-\xi_k^{n+1}|$ is often too small to control the error of $\phi$ or other variables. Fortunately, extensive numerical experiments show
$|1-\xi_k^{n+1}|^m$ performs well  when  $m$ is appropriately chosen between 0 and 1.\par
Based on these observations, we utilize the aforementioned time-stepping strategies to design adaptive algorithms for the R-IMEX-BDF$k$ scheme.
\subsection{Adaptive time-stepping strategy I}For the final time $T$ and the R-IMEX-BDF$k$ scheme, given default safety coefficients $\rho_{k}$, reference tolerances $tol_{k}$, tunable constants for adaptivity speed $r_{k}$,\
error order  constants $m_{k}$,  minimum time steps $\tau_{kmin}\ $,\  maximum time steps $\tau_{kmax}$,\
 the energy coefficient $\gamma^{*}_k$,\  and  the subscript
$k$ represents the order. The time steps of the R-IMEX-BDF$k$ scheme are updated to improve numerical accuracy by using the following formula \cite{huang2020highly}:
\begin{eqnarray*}
	A_{dp}(e,\tau)=\rho_k(\frac{tol_k}{e})^{r_k}\tau.
\end{eqnarray*}
The update algorithm is summarized as follows.
	 \begin{algorithm}[!ht]
     \caption{(R-IMEX-BDF$k$A)Adaptive time strategy}
     \label{algo1}
	 	\begin{algorithmic}
	\STATE\textbf{Given:}\hspace{0.5em}The previous time step $\tau_{n}$.
	 		\STATE\textbf{Step 1:}\hspace{0.5em}Compute $\xi_k^{n+1}$ by the Scheme 1 with time step $\tau_{n}$ ;
	 		\STATE\textbf{Step 2:}\hspace{0.5em}Calculate $e_{n}=|1-\xi_k^{n+1}|$ ;
	 		\STATE\textbf{Step 3:}\hspace{0.5em}\textbf{if} $e_{n}^{m_k}>tol_k$ ,\textbf{then}
	\STATE\hspace{5.0em} calculate  $\tau_{n}\leftarrow max\left\{ \tau_{kmin},min \left\{A_{dp}(e^{m_k}_{n},\tau_{n}),\dfrac{\tau_{kmax}}{\sqrt{1+\left|\gamma^{*}_kE^{'}(t) \right|  ^{2}}}\right\}\right\}$;
\STATE\hspace{5.0em} \textbf{goto} step 1
	 		\STATE\textbf{Step 4:}\hspace{0.5em}\textbf{else} \STATE\hspace{5.0em} update  $\tau_{n+1}\leftarrow max\left\{ \tau_{kmin},min\left\{ A_{dp}(e^{m_k}_{n},\tau_{n+1}),\dfrac{\tau_{kmax}}{\sqrt{1+\left|\gamma^{*}_kE^{'}(t) \right|^{2}}}\right\}\right\}$;
	 		\STATE\textbf{Step 5:}\hspace{0.5em}\textbf{endif}
	 	\end{algorithmic}
	 \end{algorithm}
\begin{rem}
For the case $\gamma^{\star}_k = 0$ and $m_k = 1$, the adaptive time strategy is similar to that in \cite{huang2020highly}. When $\gamma^{}_k \neq 0$ and $tol_k$ is sufficiently large, the strategy resembles those in \cite{zhang2013adaptive, hou2023implicit, cheng2018multiple}. Abundant numerical experiments demonstrate that, with appropriate parameters, our generalized time-adaptive strategy improves the efficiency of solving the CHB model.
\end{rem}

In recent works \cite{liao2022adaptive,hou2023implicit,hou2023linear}, the authors cleverly developed a series of unconditionally energy-stable schemes with rigorous analysis for some gradient flow models. The modified discrete energy remains unconditionally stable only when the ratio of adjacent time steps, $\gamma^{*} = \tau_{n+1}/\tau_{n}$, is $\leq 4.8645$. Fortunately, from Theorem \ref{chbR5}, the adjacent time steps are not restricted in this paper, and the adaptive time strategy remains unaffected by them.
As noted in \cite{jiang2024highly}, an adaptive time algorithm, IMEX-CNA, was first proposed for the CHB model. The authors use the numerical solution of the first-order IMEX-SAV method as the reference and take the relative error of the IMEX-CN solution as the error indicator, which results in high computational costs for obtaining available solutions.

As is well known, in the early stage of dynamics of some phase-field models, the energy decays rapidly due to nonlinear interactions, then more slowly until reaching a steady state \cite{zhang2013adaptive}. A similar behavior is observed in the CHB model. Typically, $\phi$ and $E_1(\phi)$ eventually reach a stationary state \cite{conti2020well}, but change rapidly in a short time interval. Therefore, a high-order time discretization scheme is recommended to capture the model's rapid evolution during this period. As the order of the BDF scheme increases, the stability region decreases.  In most cases, high-order schemes are not suitable for long-time simulations of the CHB model.
 Several tests show that the R-IMEX-BDF$k$ scheme or R-IMEX-BDF$k$A algorithm ($k \geq 3$) may not be effective for the CHB model unless small time steps are used for long simulations or large time steps for short intervals. Additionally, R-IMEX-BDF$2$A is less effective than R-IMEX-BDF$3$A in short intervals but more efficient for long-time simulations. A similar behavior was observed in \cite{huang2020highly}, where higher-order BDF schemes require more time to solve the Cahn–Hilliard equation. Based on these observations, we propose a hybrid-order adaptive time-stepping strategy for updating the time steps.
\subsection{Adaptive time-stepping strategy II}For the final time $T$, given a time node $T_c \in [0,T]$, the other parameters are the same as in Algorithm \ref{algo1}.
The algorithm is as follows.
\begin{algorithm}[!ht]
	\caption{ Adaptive time strategy}
	\label{algo2}
	\begin{algorithmic}
		\STATE \textbf{Given:} The previous time step $\tau_{n}$.
		\STATE \textbf{Step 1:} \textbf{if} $t_{n+1} \leq T_{c}$ , \textbf{then}
	\STATE\hspace{5.3em}compute $\xi^{n}$ by the R-IMEX-BDF3 with time step $\tau_{n}$ ;
		\STATE \textbf{Step 2:} \hspace{1.5em}calculate $e_{n}=|1-\xi_3^{n}|$ ;
		\STATE \textbf{Step 3:} \hspace{1.5em}\textbf{if} $e^{m_3}_{n} > tol_{3}$ , \textbf{then}
	\STATE \hspace{6.5em}	recalculate  $\tau_{n} \leftarrow \max\left\{\tau_{3min}, \min\left\{ A_{dp}(e_{n}^{m_3}, \tau_{n}), \frac{\tau_{3max}}{\sqrt{1+\left|\gamma^{*}_3E' (t) \right|  ^{2}}}\right\}\right\}$;
	\STATE \hspace{6.5em} \textbf{goto} Step 1
		\STATE \textbf{Step 4:} \hspace{1.5em}\textbf{else}
		\STATE \hspace{6.5em} update  $\tau_{n+1} \leftarrow \max\left\{ \tau_{3min}, \min\left\{ A_{dp}(e^{m_3}_{n}, \tau_{n}),\frac{\tau_{3max}}{\sqrt{1+\left|\gamma^{*}_3E' (t) \right|  ^{2}}}\right\}\right\}$;
		\STATE \textbf{Step 5:} \hspace{1.5em}\textbf{endif}
		\STATE \textbf{Step 6:} \textbf{else}
	\STATE\hspace{5.3em} compute $\xi^{n}$ by the R-IMEX-BDF2 with  $\tau_{n}$ ;
		\STATE \textbf{Step 7:} \hspace{1.5em} calculate $e_{n}=|1-\xi_2^{n}|$;				
		\STATE \textbf{Step 8:} \hspace{1.5em} \textbf{if} $e^{m_2}_{n} > tol_{2}$ , \textbf{then}
	\STATE \hspace{6.5em} recalculate $\tau_{n} \leftarrow \max\left\{\tau_{2min}, \min\left\{ A_{dp}(e_{n}^{m_2}, \tau_{n}), \frac{\tau_{2max}}{\sqrt{1+\left|\gamma^{*}_2E' (t) \right|  ^{2}}}\right\}\right\}$;
	\STATE \hspace{6.5em} \textbf{goto} Step 6
		\STATE \textbf{Step 9:} \hspace{1.5em} \textbf{else}
	\STATE \hspace{6.5em} update  $\tau_{n+1} \leftarrow \max\left\{ \tau_{2min}, \min\left\{ A_{dp}(e_{n}^{m_2}, \tau_{n}), \frac{\tau_{2max}}{\sqrt{1+\left|\gamma^{*}_2E' (t) \right|  ^{2}}}\right\}\right\}$;
		\STATE\textbf{Step 10:} \hspace{1em} \textbf{endif}
		\STATE \textbf{Step 11:} \textbf{endif}
	\end{algorithmic}
\end{algorithm}
\begin{rem}
For $T_c = 0$, Algorithm \ref{algo2} corresponds to the R-IMEX-BDF2A method on $[0,\ T]$. For $T_c = T$, Algorithm \ref{algo2} becomes the R-IMEX-BDF3A method on $[0,\ T]$. For $T_c \in (0, T)$, Algorithm \ref{algo2} is a hybrid-order adaptive time-stepping algorithm on $[0, T]$. By choosing an appropriate $T_c$, Algorithm \ref{algo2} applies R-IMEX-BDF3A for high accuracy with large time steps in the short interval $[0,\ T_c]$, while R-IMEX-BDF2A is used to maintain relatively large steps for the longer simulation in $(T_c,\ T]$. Several experiments show that Algorithm \ref{algo2} effectively combines R-IMEX-BDF2A and R-IMEX-BDF3A, achieving high efficiency for the CHB model. It is worth noting that the Algorithm \ref{algo2} does not alter the fact that the modified energy $r^n$ remains unconditionally decaying.
\end{rem}

\section{Numerical experiments}
\setcounter{equation}{0}
To demonstrate the performance of the proposed numerical method, we present several numerical experiments in this section. For convenience, we assume periodic boundary conditions on the computational domain $\Omega = [0, 2\pi]^2$. We use $256^2$ Fourier modes, which are sufficiently large to make the spatial discretization errors negligible compared to the time discretization errors.  For simplicity, we set $C_0=0$.
In fact, similar results can be obtained for general values of
$C_0$ in the following numerical experiments.
All numerical experiments were carried out on a computer with the following
characteristics: Intel(R) Core(TM)i9-12900k CPU@3.20GHz, RAM 128 GB.
\begin{example}(Accuracy test)
In this first example, we test the convergence rates of the R-IMEX-BDF$k$ method with  fixed time steps for the CHB model by using the  following exact solution:
\begin{eqnarray}\label{exp1}
\begin{array}{ll}
\phi(x,y,t)=cos(x)sin(y)cos(t),\\
u_{1}(x,y,t)=sin(x)sin(y)sin(t),\\
u_{2}(x,y,t)=cos(x)cos(y)sin(t),\\
p(x,y,t)=cos(x)sin(y)sin(t).\\
\end{array}
\end{eqnarray}
We set the following parameters:
\begin{equation*}
\varepsilon=1 , \gamma=2 , M=\nu=\eta=1 , S=0.
\end{equation*}
The $L^2$ errors and convergence order of the R-IMEX-BDF$k$ ($k = 1, 2, 3, 4$) are shown in Figure \ref{fig1}. All $L^2$ errors are computed by comparing the numerical solution with the exact solution \eqref{exp1} at the final time $T = 1$. As expected, the results confirm the temporal convergence order for $\phi$, $\mathbf{u}$, $p_x$, and $p_y$ of \textbf{Scheme 1}.
\begin{figure}[!ht]
\centering
\subfigure[ R-IMEX-BDF1 scheme.]{
\includegraphics[width=7.3cm]{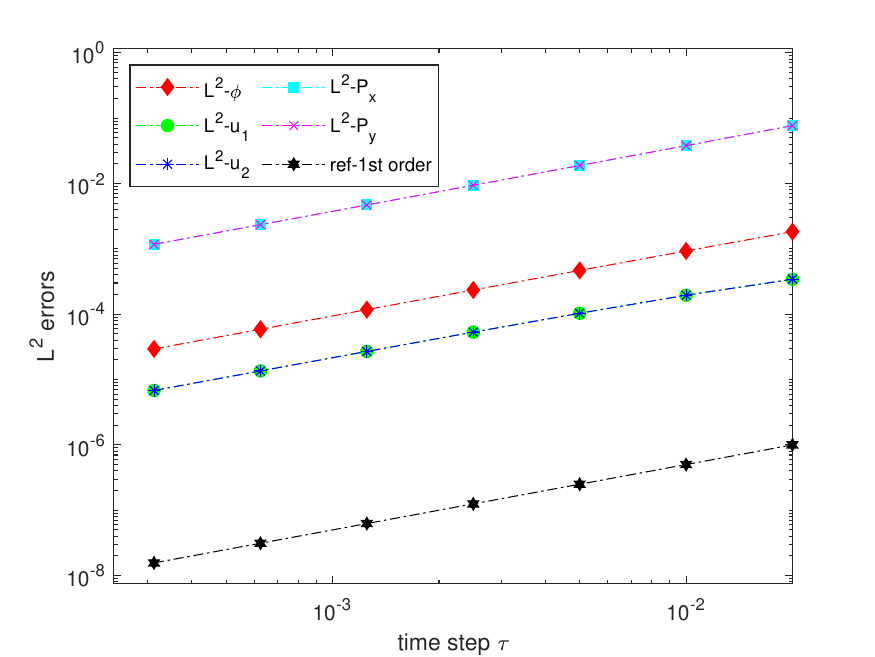}
}
\hspace{-10mm}
\subfigure[R-IMEX-BDF2 scheme.]{
\includegraphics[width=7.3cm]{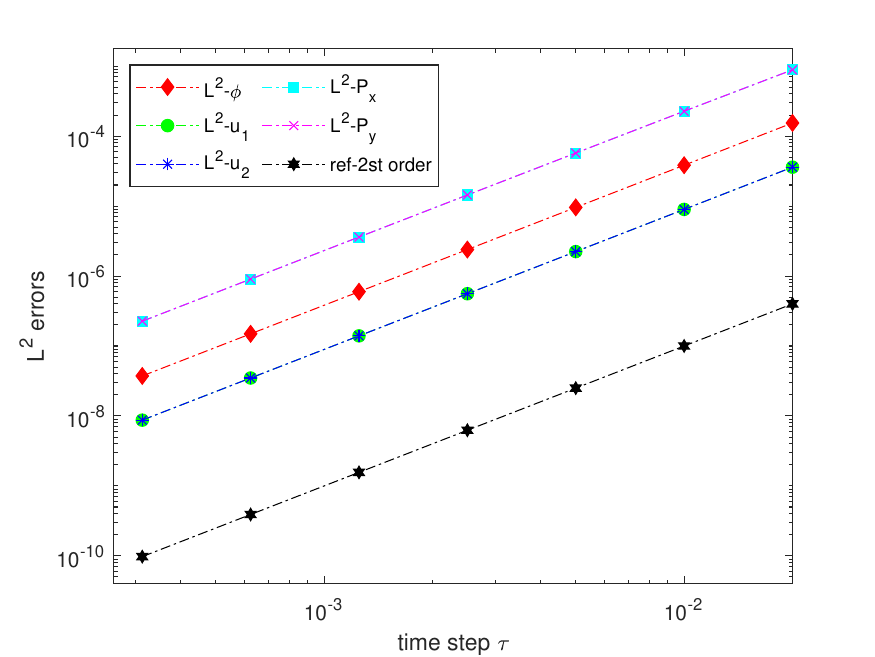}}
\quad
\subfigure[R-IMEX-BDF3 scheme.]{
\includegraphics[width=7.3cm]{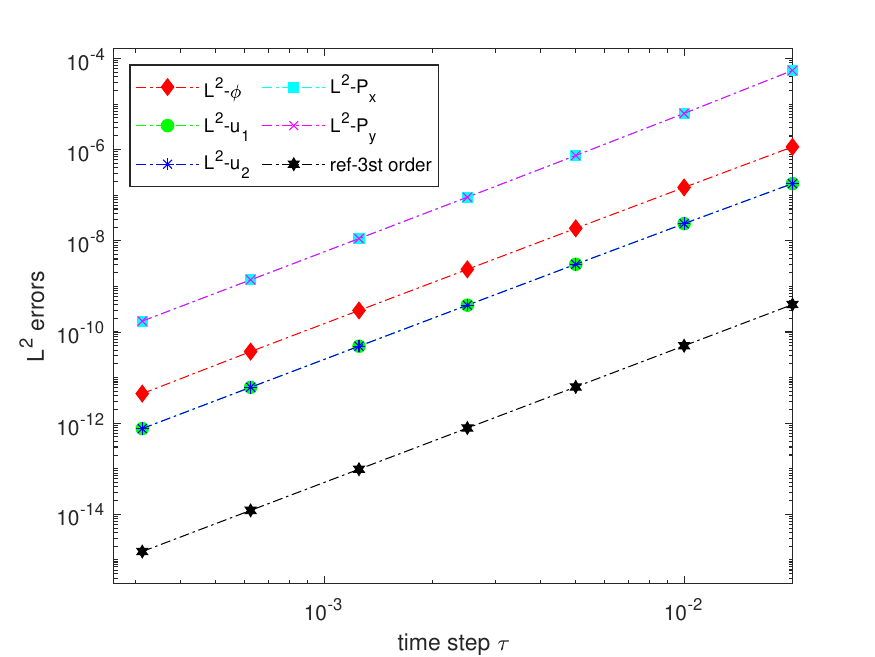}
}
\hspace{-10mm}
\subfigure[R-IMEX-BDF4 scheme.]{
\includegraphics[width=7.3cm]{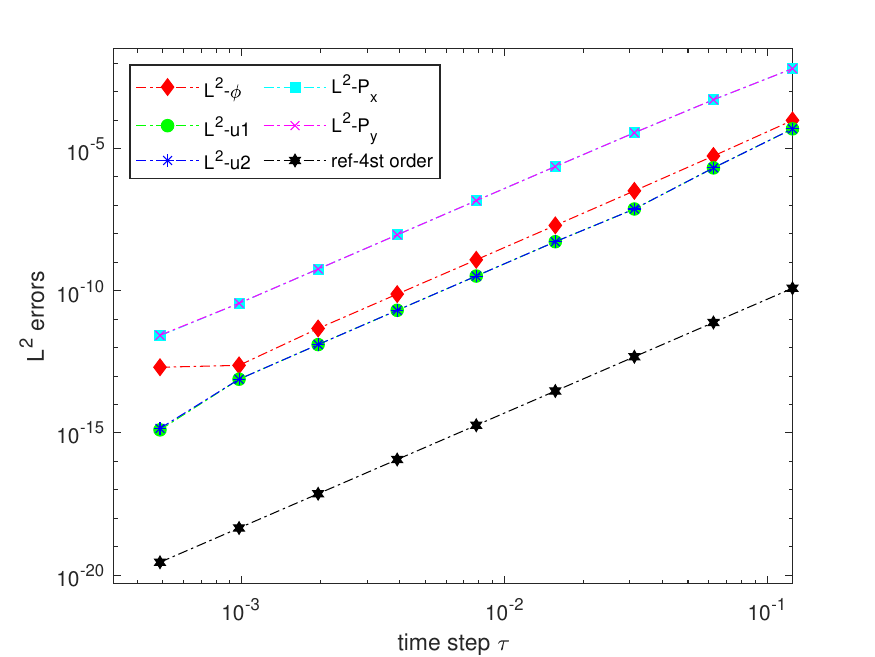}
}
\caption{(Example \ref{exp1}) Temporal convergence rates: $L^2$ errors of phase field function $\phi$, the velocity ($x$-component of $u_{1}$, $y$-component of $u_{2}$), and the derivative of pressure ($p_{x}$ and $p_{y}$)  as functions of the time step size $\tau$ for \textbf{Scheme 1}.}
\label{fig1}
\end{figure}

\end{example}
\begin{example}\label{exp2} (The effect of relaxation and  coarsening process)
In the second example, we simulate the classical coarsening process of a two-phase homogeneous mixture.  The process occurs when a homogeneous system is disturbed by an external perturbation. This helps verify the effectiveness of relaxation in the R-IMEX-BDF$k$ scheme for the CHB model. We initialize the CHB system with the following conditions:
\begin{eqnarray}
\label{im111}
\begin{array}{ll}
\phi(x,y,0)=-0.5-0.001*(2*Rand(x,y)-1),
\end{array}
\end{eqnarray}
where the following parameters are used in this test:
\begin{equation*}
	\varepsilon=0.05 , \gamma=4 , M=\nu=\eta=1,\ S=1.
\end{equation*}

We use the results from the semi-implicit/BDF2 scheme with \( \tau_n = 1e-6 \) as the reference solution. The \( L^2 \)-norm errors of the R-IMEX-BDF2 scheme and the second-order IMEX-BDF2 scheme at \( T = 100 \) with various time steps are shown in Table \ref{tab1}. We observe that the R-IMEX-BDF2 scheme reduces the \( L^2 \)-norm error compared to IMEX-BDF2, improving accuracy, especially for larger time steps.

Figure \ref{fig61}  presents the time histories of the modified energy $r(t)$ using large time step sizes $\tau_n=1$ and $\tau_n=4$. While such large time steps may compromise accuracy, the modified energy consistently decays and remains positive throughout the long simulations, indicating the schemes' robustness. It illustrates that the proposed R-IMEX-BDF$k$ ($k=2,3$) schemes with fixed time steps keep unconditional energy stable.

In Figure \ref{fig2}(a), the energy history curves over a long period are plotted, showing that the curves for \( \tau_n = 1e-3 \) overlap and consistently exhibit energy decay. Figures \ref{fig2}(b) and (c) display the error of the energy  and \( \xi \) , respectively, where it is evident that the R-IMEX-BDF2 scheme improves the accuracy of energy and significantly reduces the error of \( \xi \). These observations align with the relaxation behavior seen in the Allen-Cahn equation in \cite{zhang2022generalized}.
\begin{table}[!ht]
	\centering
    	\caption{(Example \ref{exp2}) A comparison of the \( L^2 \)-norm errors of $\phi$ between the R-IMEX-BDF2 and IMEX-BDF2 schemes  at \( T = 100 \) with different time steps.
    }
	\begin{tabular}{p{2.7cm} p{2.7cm} p{2.7cm} p{2.7cm} p{2.7cm}}
		\bottomrule
		 & $1e-3$ & $1.5e-3$ & $1.8e-3$ & $2e-3$\\
		\midrule
		IMEX-BDF2 & 0.0934 & 0.2463  & 0.4178 & 2.0342\\
		R-IMEX-BDF2 & 0.0922 & 0.2359 & 0.3712 & 0.4870\\
		\bottomrule
	\end{tabular}
		\label{tab1}
\end{table}
\begin{figure}[!ht]
	\centering
	\subfigure[R-IMEX-BDF2]{
	\includegraphics[width=7.3cm]{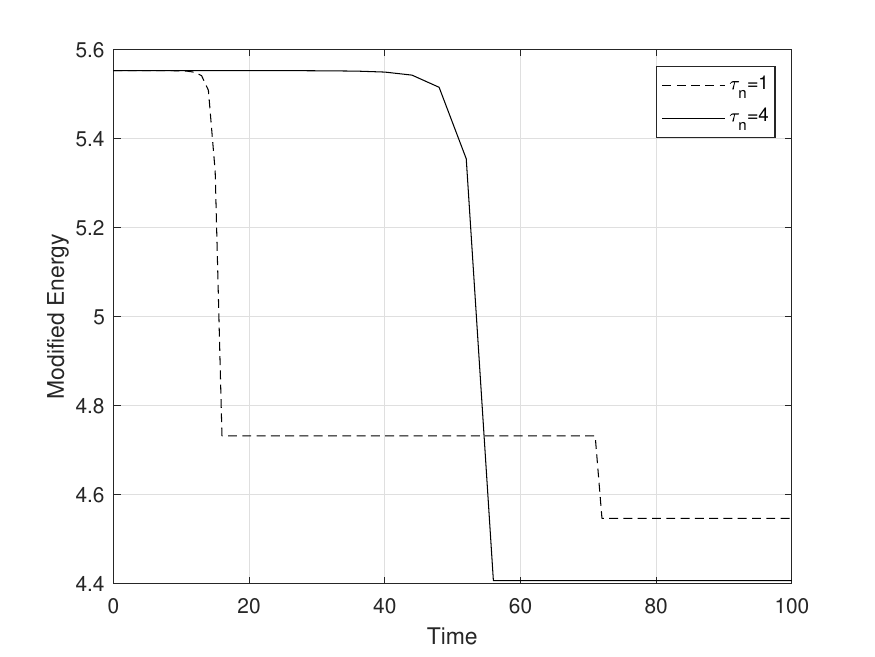}
}
\hspace{-10mm}
\subfigure[R-IMEX-BDF3]{
	\includegraphics[width=7.3cm]{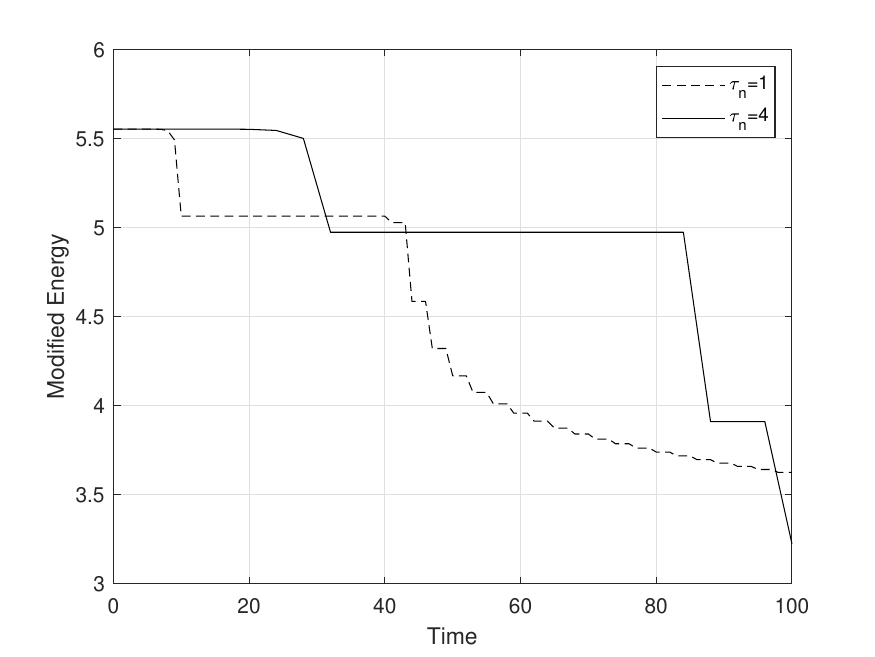}
}
\caption{(Example \ref{exp2})Time histories of the modified energy computed by (a) the R-IMEX-BDF2 scheme, (b) the R-IMEX-BDF3 scheme,using large time step sizes $\tau_n$=1,4.}
	\quad
\label{fig61}
\end{figure}

\begin{figure}[!ht]
\centering
\subfigure{
\includegraphics[width=2.1in]{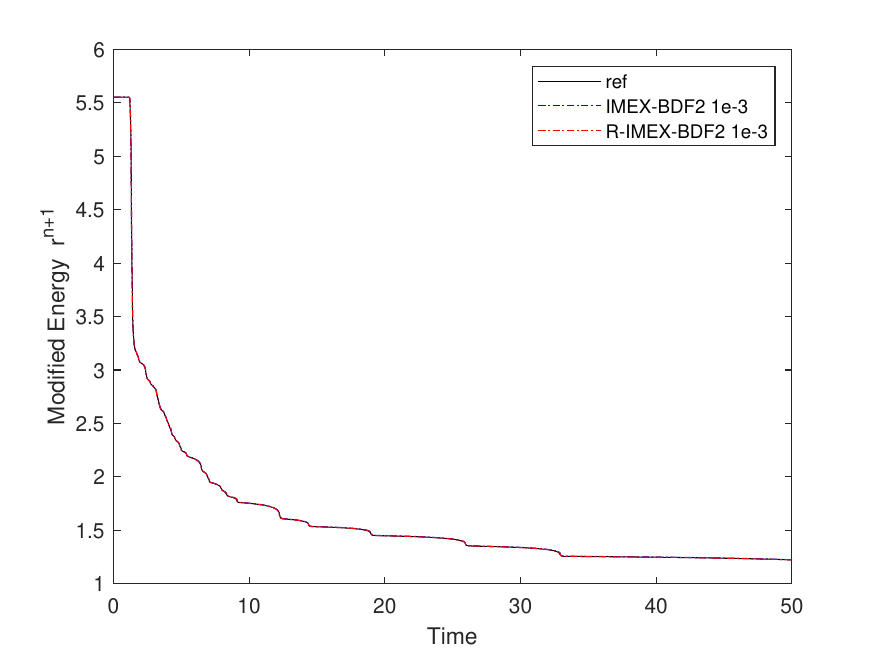}
}
\hspace{-8mm}
\subfigure{
\includegraphics[width=2.1in]{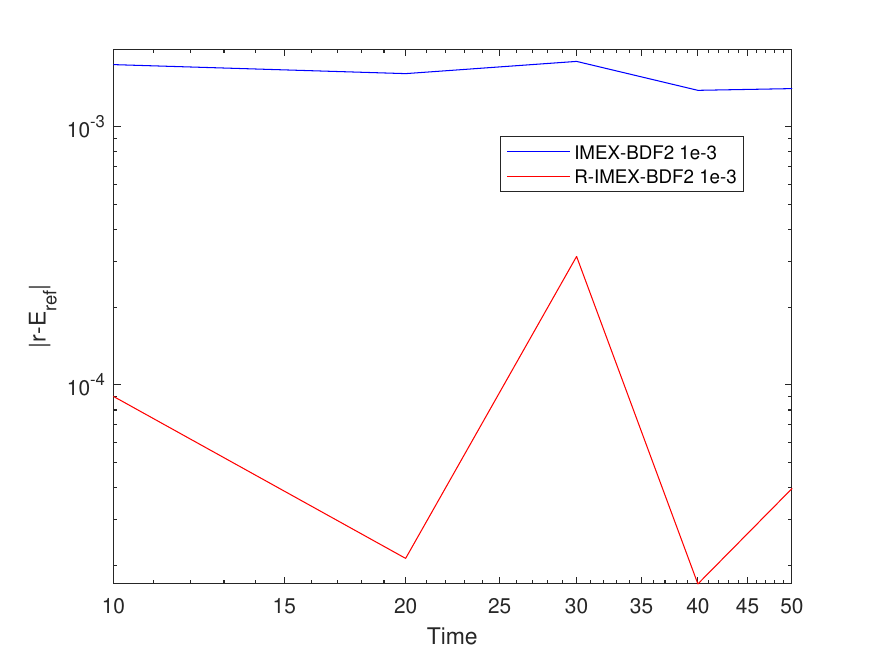}
}
\subfigure{
\includegraphics[width=2.1in]{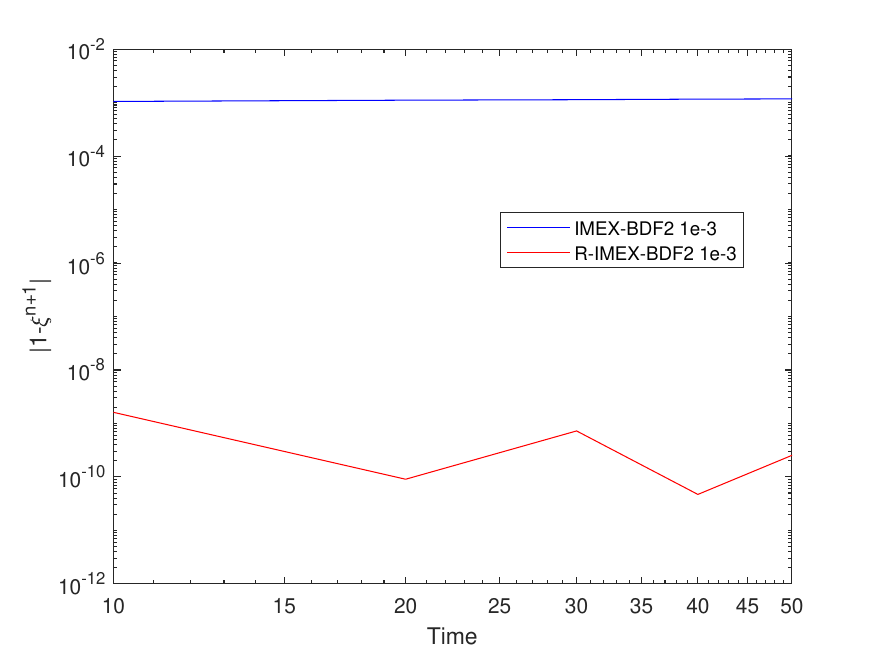}
}
\caption{(Example \ref{exp2})Time evolutions of modified energy (left), error (middle), and  $e_{n+1}$ (right).}
\label{fig2}
\end{figure}

In addition, we simulate the spinodal decomposition of a homogeneous mixture into two coexisting phases governed by the CHB system as a further test of the algorithms developed here. Figure \ref{fig6} shows the results of the coarsening process for initial concentration \( \phi_0 = -0.5 \) at times \( T = 1.5, 4, 8, 20, 80,150, 200, 500 \). At \( T = 1.5 \), the two-phase mixture forms many small balls scattered throughout the domain. As time progresses, the number of these balls decreases while their size increases, consistent with the observations in \cite{jiang2024highly}.
	\begin{figure}[!ht]
		\centering
		\subfigure[$T=1.8$]{
			\begin{minipage}[t]{0.2\linewidth}
				\includegraphics[width=1.5in]{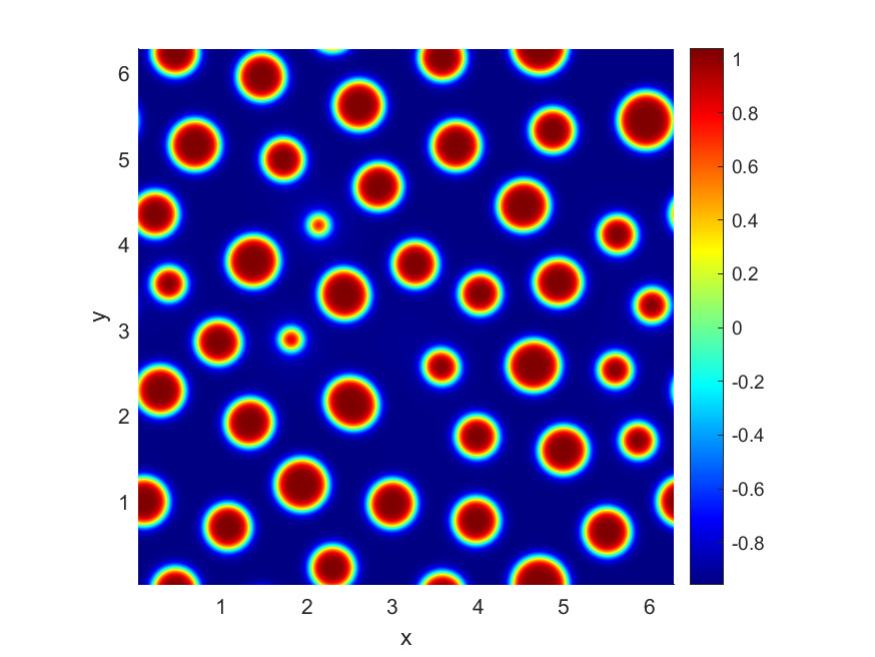}
			\end{minipage}%
		}%
		\subfigure[$T=4$]{
			\begin{minipage}[t]{0.2\linewidth}
				\includegraphics[width=1.5in]{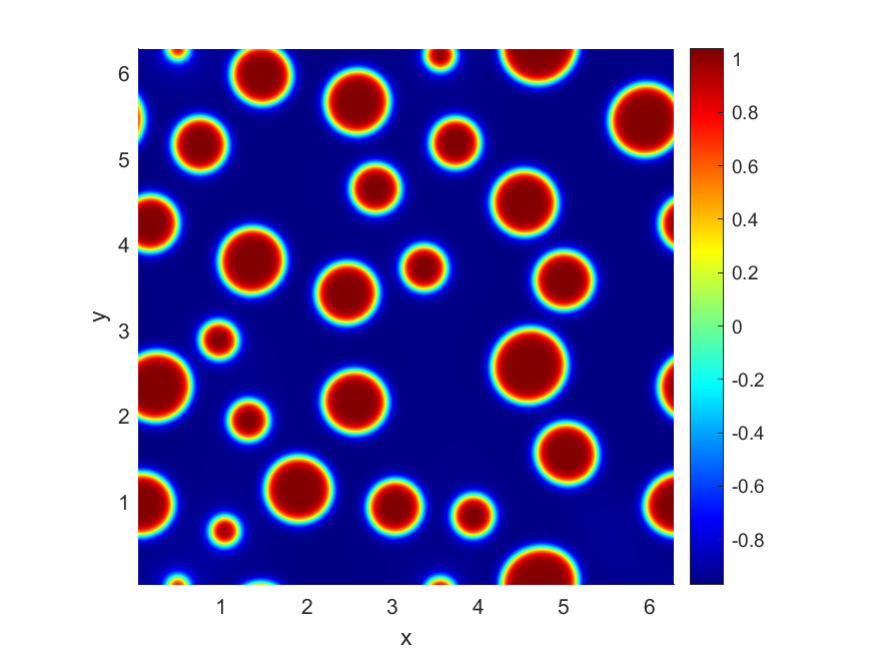}
			\end{minipage}%
		}%
		\subfigure[$T=8$]{
			\begin{minipage}[t]{0.2\linewidth}
				\includegraphics[width=1.5in]{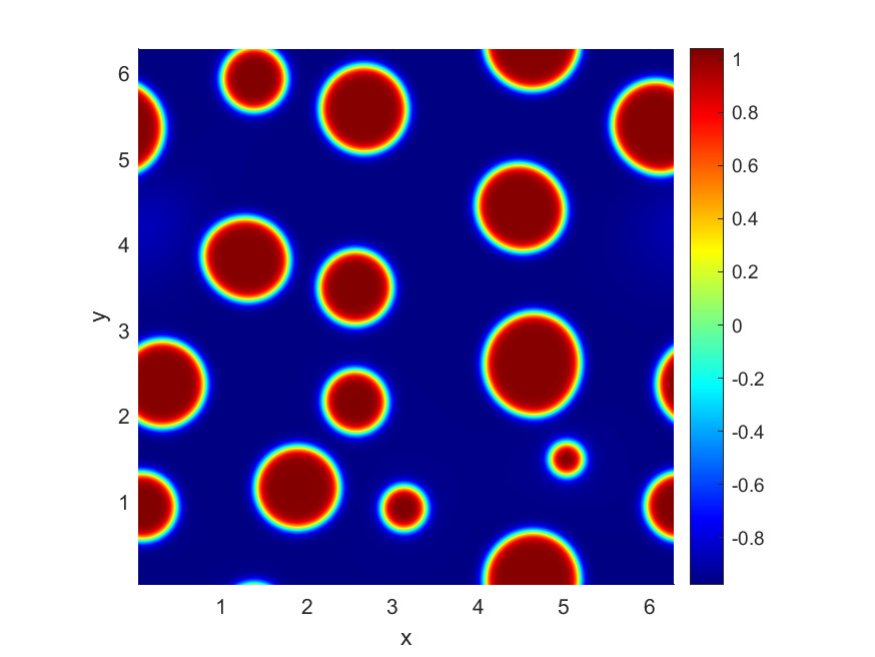}
			\end{minipage}%
		}%
		\subfigure[$T=15$]{
			\begin{minipage}[t]{0.2\linewidth}
				\includegraphics[width=1.5in]{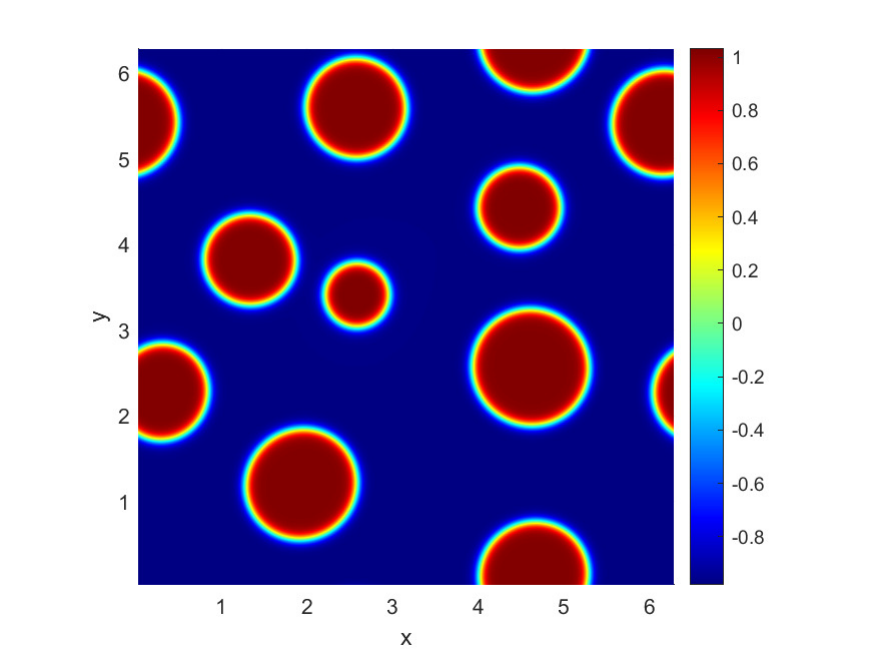}
			\end{minipage}%
		}%
		\quad
		\subfigure[$T=80$]{
			\begin{minipage}[t]{0.2\linewidth}
				\includegraphics[width=1.5in]{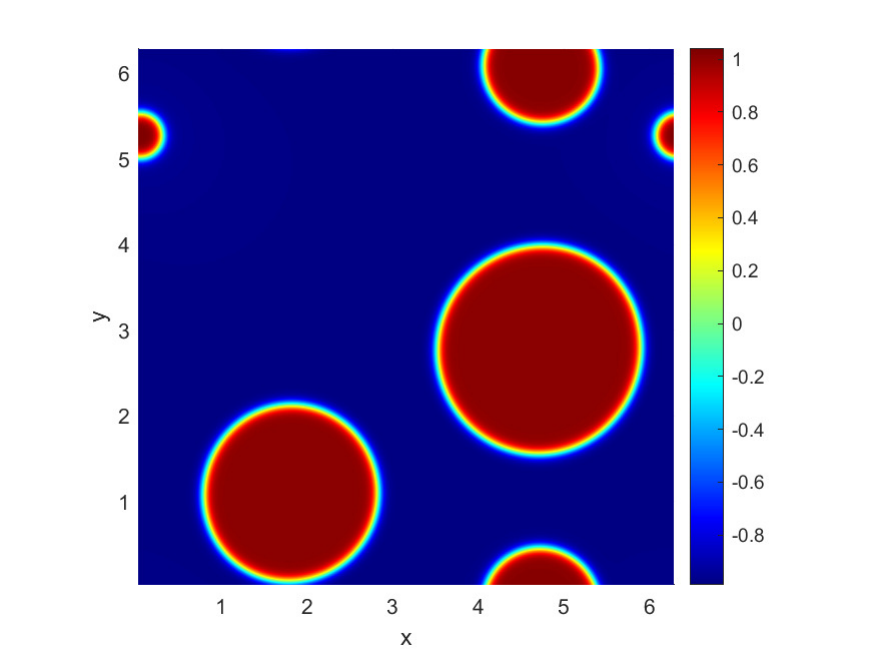}
			\end{minipage}%
		}%
		\subfigure[$T=150$]{
			\begin{minipage}[t]{0.2\linewidth}
				\includegraphics[width=1.5in]{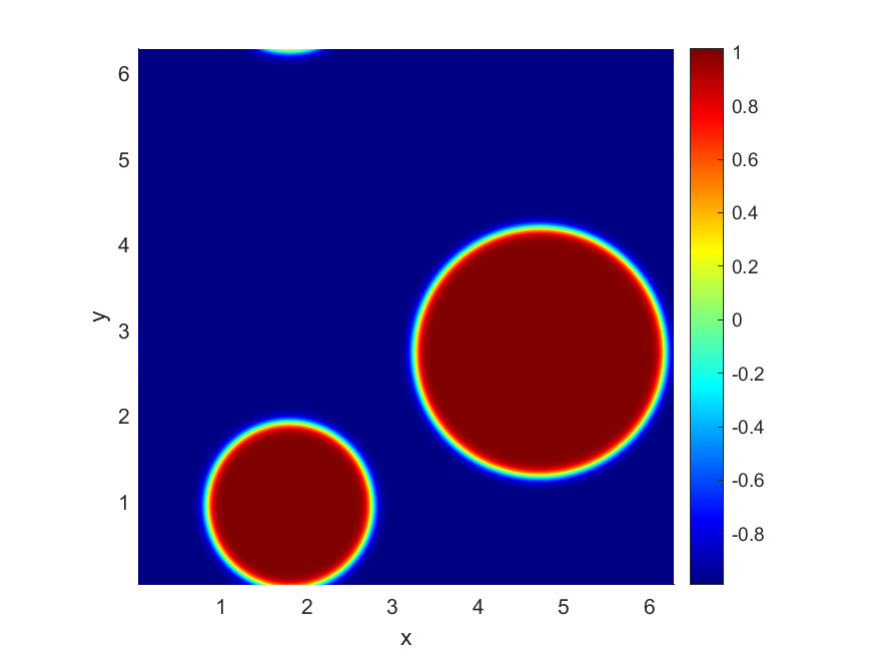}
			\end{minipage}%
		}%
		\subfigure[$T=200$]{
			\begin{minipage}[t]{0.2\linewidth}
				\includegraphics[width=1.5in]{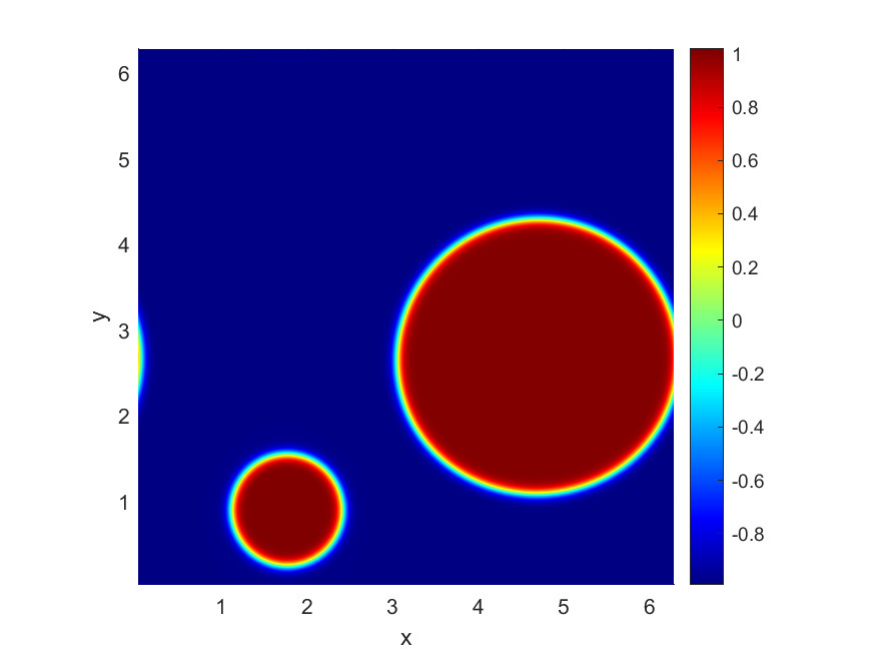}
			\end{minipage}%
		}%
		\subfigure[$T=500$]{
			\begin{minipage}[t]{0.2\linewidth}
				\includegraphics[width=1.5in]{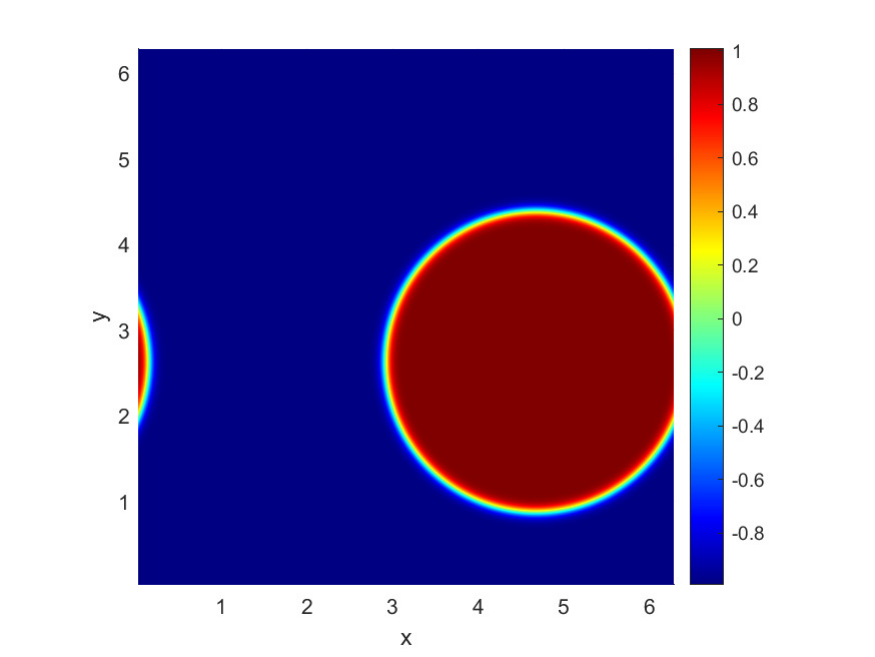}
			\end{minipage}%
		}%
		\centering
		\caption{(Example \ref{exp2})Time evolution of the coarsening process with $\epsilon=0.05$ by using the R-IMEX-BDF2 scheme with \( \tau_n = 0.001 \).}
		\label{fig6}
	\end{figure}	
\end{example}
\begin{example}(The  adaptive time steps strategy) \label{exp3}
In this example, we assess the impact of the proposed adaptive methods on the CHB model, focusing on spinodal decomposition as described in Example \ref{exp2}. Extensive experiments indicate that numerical results exhibit similar behaviors when \( S \) and \( \gamma^{*} \) are not large.   Consequently,  for the following adaptive algorithms herein, we set \( S = 1 \) and \( \gamma^{*} = 1 \).  Similar to the Example \ref{exp2}, we use the results from the R-IMEX-BDF2 scheme with $\tau_n= 1e-6$ as the reference solution.
To present the experimental results concisely, we use six representative parameter sets as follows:
\begin{eqnarray}
	&&\rho=0.75,tol=1e-3,r=0.50,m=0.52,\tau_{max}=3e-3;\label{par6}	\\
	&&\rho=0.70,tol=1e-3,r=0.70,m=0.70,\tau_{max}=4e-3;	\label{par7}\\
  &&\rho=0.68,tol=1e-3,r=0.72,m=0.75,\tau_{max}=4e-3;\label{par3}\\
	&&\rho=0.90,tol=1e-4,r=0.95,m=0.85,\tau_{max}=1.5e-3.\label{par2}\\
    &&\rho=0.50,tol=1e-3,r=0.80,m=0.68,\tau_{max}=1e-1.\label{par10}\\
    &&\rho=0.60,tol=8e-4,r=0.60,m=0.70,\tau_{max}=1e-2.\label{par11}
\end{eqnarray}
\par
In Table \ref{tab4}, we present the $L^{2}$-norm errors of $\phi$ and CPU times for R-IMEX-BDF2A and R-IMEX-BDF3A at $T=1.2$, $T=10$, and $T=50$. The numerical parameters are as specified in (\ref{par6}) and (\ref{par7}). These results highlight the numerical characteristics of the two algorithms.

At $T=1.2$, R-IMEX-BDF3A improves accuracy and reduces CPU time by a factor of eight compared to R-IMEX-BDF2 scheme with fixed time-step $\tau_n=5e-4$.  Relative to R-IMEX-BDF2A, R-IMEX-BDF3A achieves nearly a 100-fold accuracy improvement at an additional cost of 3.8 seconds. In this phase, R-IMEX-BDF3A outperforms both R-IMEX-BDF2A and fixed time-step R-IMEX-BDF2. Figures \ref{fig3}(a)-(c) show that energy decays rapidly over $(0,T)$. However, as the phase separation over time (e.g., at $T=10$ or $T=50$), the efficiency of R-IMEX-BDF3A deteriorates. At $T=50$, the error and CPU time increase substantially.
In contrast, R-IMEX-BDF2A demonstrates stable performance. When extending the simulation from $T=10$ to $T=50$, its CPU time increases by less than 100 seconds. Note that, for the same problem, higher-order BDF schemes usually have smaller stability regions. For the CHB equation, even  R-IMEX-BDF$k$ ($k\geq3$) schemes with constant time steps yield meaningful solutions only when $\tau_n\leq 2e-4$ at $T=50$ and this is not suitable for long-time simulations.

The experiment indicates that R-IMEX-BDF3A is advantageous in the early stage of phase evolution but is less suitable for long-time simulations. Conversely, R-IMEX-BDF2A, despite lower accuracy, offers stronger stability for long-term evolution. Moreover, under equivalent conditions, the R-IMEX-BDF$k$A($k\geq4$) algorithm fails to produce meaningful numerical solutions at $T=0.1$.
	\begin{table}[!ht]
		\centering
\caption{(Example \ref{exp3})The $L^{2}$-norm error of $\phi$ and CPU time for R-IMEX-BDF$k$A($k=2,3$) over time.}
		\begin{tabular}{p{6cm} p{3cm} p{2cm} p{2cm}}
			\bottomrule
			Scheme/Algorithm & $L^{2}$-norm error & CPU time &  Time $T$\\
			\midrule
			R-IMEX-BDF3A with (\ref{par6})& 7.473e-4 & 8.4s & \ \ \ 1.2 \\
			R-IMEX-BDF2A with (\ref{par7})& 6.179e-2 & 4.6s & \ \ \ 1.2\\
			R-IMEX-BDF2 with $\tau_n=5e-4$ & 1.092e-3 & 73.8s & \ \ \ 1.2\\
			\bottomrule
  R-IMEX-BDF3A with (\ref{par6})& 4.797e-3 & 368.0s & \ \ \ 10\\
     R-IMEX-BDF2A with (\ref{par7}) & 1.716e-1 &  59.2s & \ \ \ 10 \\
R-IMEX-BDF2 with $\tau_n=5e-4$& 2.936e-3 & 394.7s & \ \ \ 10\\
  \bottomrule
R-IMEX-BDF3A with (\ref{par6}) & 6.820e-2 & 2115.8s & \ \ \ 50 \\
R-IMEX-BDF2A with (\ref{par7}) & 6.433e-1 & 258.9s & \ \ \ 50 \\

R-IMEX-BDF2 with $\tau_n=5e-4$ & 7.811e-3 & 1855.9s & \ \ \ 50 \\
  \bottomrule
\end{tabular}

		\label{tab4}
\end{table}
\par
Algorithm \ref{algo2} combines the advantages of R-IMEX-BDF2A and R-IMEX-BDF3A. It reduces CPU time while maintaining high accuracy. In the early stage \([0, T_c]\), R-IMEX-BDF3A is employed to achieve high accuracy with larger time steps. Subsequently, R-IMEX-BDF2A is utilized to sustain larger steps for the extended simulation period \((T_c, T]\). For this study, we set \(T_c = 1.2\) to evaluate the performance of this combined approach.
\par
In Table \ref{tab5}, at \(T=10\), Algorithm \ref{algo2} with parameters (\ref{par6}) and (\ref{par7}) improves accuracy by at least eight times compared to R-IMEX-BDF2A with parameter (\ref{par7}). The CPU time increases by only 10 seconds. Compared to a fixed time step of \(\tau_n=2e-3\), it achieves higher accuracy with lower cost.
At \(T=50\), the advantage becomes more significant. Compared to R-IMEX-BDF2A, Algorithm \ref{algo2} improves accuracy by at least nine times while reducing CPU time by 36 seconds. Compared to a fixed time step of \(\tau_n=2e-3\), it achieves more than twice the accuracy and cuts CPU time by over threefold.
\begin{table}[!ht]
		\centering
\caption{(Example \ref{exp3})The $L^{2}$-norm error of $\phi$ and CPU time for Algorithm \ref{algo2} and R-IMEX-BDF2A at $T=10$ and $T=50$.}
		\begin{tabular}{p{6cm} p{3cm} p{2cm} p{2cm}}
			\bottomrule
			Scheme/Algorithm & $L^{2}$-norm error & CPU time &  Time $T$\\
			\midrule
  Algorithm \ref{algo2} with (\ref{par6})(\ref{par7}) & 2.293e-2 & 70.2s & \ \ \ 10 \\
R-IMEX-BDF2A with (\ref{par7}) & 1.716e-1 &  59.2s & \ \ \ 10 \\
R-IMEX-BDF2 with $\tau_n=1e-3$& 1.196e-2 & 199.7s & \ \ \ 10\\
R-IMEX-BDF2 with $\tau_n=2e-3$& 5.523e-2 & 100.2s & \ \ \ 10\\
  \bottomrule
Algorithm \ref{algo2} with (\ref{par6})(\ref{par7}) & 7.131e-2 & 222.5s & \ \ \ 50 \\
R-IMEX-BDF2A with (\ref{par7}) & 6.433e-1 & 258.9s & \ \ \ 50 \\
R-IMEX-BDF2 with $\tau_n=1e-3$ & 3.275e-2 & 1184.9 & \ \ \ 50 \\
R-IMEX-BDF2 with $\tau_n=2e-3$ & 1.530e-1 & 722.5s & \ \ \ 50 \\
  \bottomrule
\end{tabular}
		\label{tab5}
\end{table}
\par
In Table \ref{tab6}, at \( T = 10 \), both R-IMEX-BDF2A and Algorithm \ref{algo2} outperform the IMEX-CNA method \cite{jiang2024highly} for the CHB equation. With parameters (\ref{par2}), both R-IMEX-BDF2A and Algorithm \ref{algo2} surpass IMEX-CNA in accuracy and CPU time. Notably, IMEX-CNA fails to produce valid numerical solutions with parameters (\ref{par3}), whereas both R-IMEX-BDF2A and Algorithm \ref{algo2} succeed. As shown in Table \ref{tab6}, Algorithm \ref{algo2} achieves  higher accuracy than R-IMEX-BDF2A under both parameters (\ref{par2}) and (\ref{par3}) with less CPU time. Furthermore, R-IMEX-BDF2A results in large errors and meaningless solutions with parameters (\ref{par10}) and (\ref{par11}), while Algorithm \ref{algo2} provides numerical solutions with acceptable errors and less computational cost.

In Table \ref{tab6},  Algorithm \ref{algo2} demonstrates superior performance compared to R-IMEX-BDF2A with fixed time steps. Compared to R-IMEX-BDF2A with fixed time steps \(\tau_n=2 \times 10^{-3}\), Algorithm \ref{algo2}  with parameters (\ref{par6}) and (\ref{par10})  achieves similar error while reducing CPU time by more than twofold. Compared to R-IMEX-BDF2A with fixed time steps  \(\tau_n = 1 \times 10^{-3}\), Algorithm \ref{algo2} with parameters (\ref{par6}) and (\ref{par11}) maintains comparable accuracy and cuts CPU time by over fourfold. Additionally, compared to R-IMEX-BDF2A with \(\tau_n = 5 \times 10^{-4}\), Algorithm \ref{algo2} with parameters (\ref{par6}) and (\ref{par2}) improves accuracy by more than twofold and reduces CPU time by over twofold. These results indicate that Algorithm \ref{algo2} consistently outperforms R-IMEX-BDF2A with fixed time steps in both accuracy and computational efficiency for the CHB equation under various parameter settings.

\begin{table}[!ht]
		\centering
\caption{(Example \ref{exp3})The $L^{2}$-norm error of $\phi$ and CPU time for some adaptive algorithms at $T=10$.}
		\begin{tabular}{p{6cm} p{3cm} p{2cm} p{2cm}}
			\bottomrule
			Scheme/Algorithm & $L^{2}$-norm error & CPU time &  Time $T$\\
			\midrule
  Algorithm \ref{algo2} with (\ref{par6})(\ref{par3}) & 1.114e-2 & 103.2s & \ \ \ 10 \\
    Algorithm \ref{algo2} with (\ref{par6})(\ref{par2}) & 1.214e-3 & 156.9s & \ \ \ 10 \\
    Algorithm \ref{algo2} with (\ref{par6})(\ref{par10}) & 5.998e-2 & 38.8s & \ \ \ 10 \\
    Algorithm \ref{algo2} with (\ref{par6})(\ref{par11}) & 1.333e-2 & 47.2s & \ \ \ 10 \\
     R-IMEX-BDF2A with (\ref{par3}) & 8.537e-2 & 105.7s & \ \ \ 10 \\
R-IMEX-BDF2A with (\ref{par2}) & 3.319e-2 & 164.9s & \ \ \ 10 \\
IMEX-CNA with (\ref{par3}) & \--\-- & \--\-- & \ \ \ 10\\
IMEX-CNA with (\ref{par2}) & 4.077e-2 & 348.3s & \ \ \ 10\\
R-IMEX-BDF2 with $\tau_n=5e-4$& 2.936e-3 & 394.7s & \ \ \ 10\\
R-IMEX-BDF2 with $\tau_n=1e-3$& 1.196e-2 & 199.7s & \ \ \ 10\\
R-IMEX-BDF2 with $\tau_n=2e-3$& 5.523e-2 & 100.2s & \ \ \ 10\\
  \bottomrule
\end{tabular}
		\label{tab6}
\end{table}
\par
In Figures \ref{fig3}(a)-(c), we present the modified energy $r^n$ over the interval $(0,500)$ for Algorithm \ref{algo2}, R-IMEX-BDF2A, and R-IMEX-BDF2, respectively.  Despite employing large adaptive time steps, the modified energy \( r^n \) consistently decreases, consistent with theoretical predictions.  Snapshots (d)-(g) at $T=500$ show very similar phase interfaces,  further demonstrating the validity of our numerical solution. Table \ref{tab7} lists the CPU times required by different algorithms at $T=500$. For long-time simulations,  Algorithm \ref{algo2} with  parameters (\ref{par6}) and (\ref{par10}) requires much less time than R-IMEX-BDF2 with $\tau_n=1e-3$.  Additionally, R-IMEX-BDF2A with parameters (\ref{par3}) reduces CPU time by approximately one-quarter. Notably, as shown in Table \ref{tab1}, when using constant time steps \(\tau_n = 2 \times 10^{-3}\), R-IMEX-BDF2 produces large errors at \(T=100\). R-IMEX-BDF2 with fixed time steps \(\tau_n \geq 2 \times 10^{-3}\) will not be suitable for long-time simulations in this case.

\begin{figure}[!ht]
\centering

\subfigure[Algorithm \ref{algo2} with (\ref{par6})(\ref{par10}).]{
\includegraphics[width=2.1in]{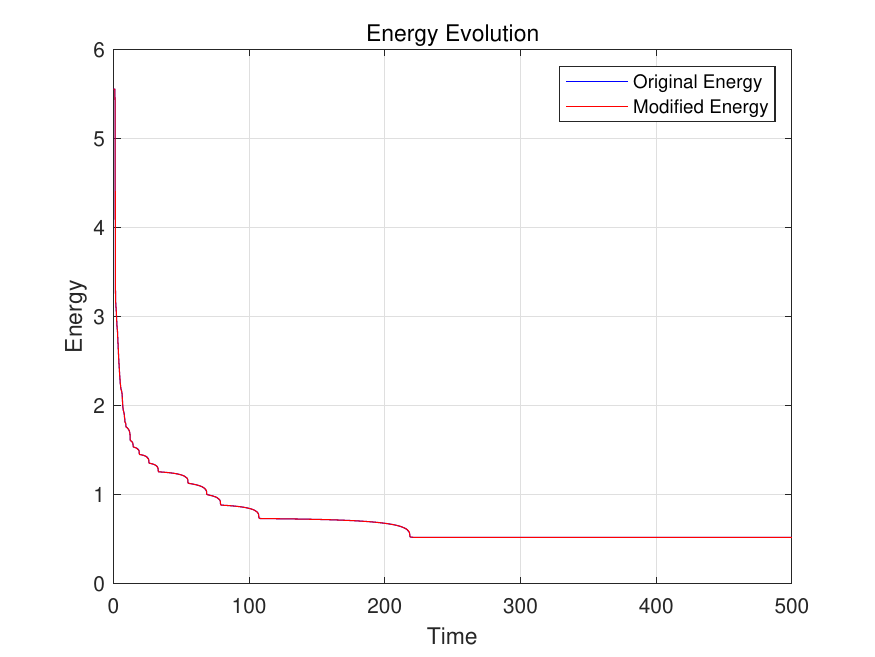}
}
\hspace{-8mm}
\subfigure[R-IMEX-BDF2A with (\ref{par3}).]{
\includegraphics[width=2.1in]{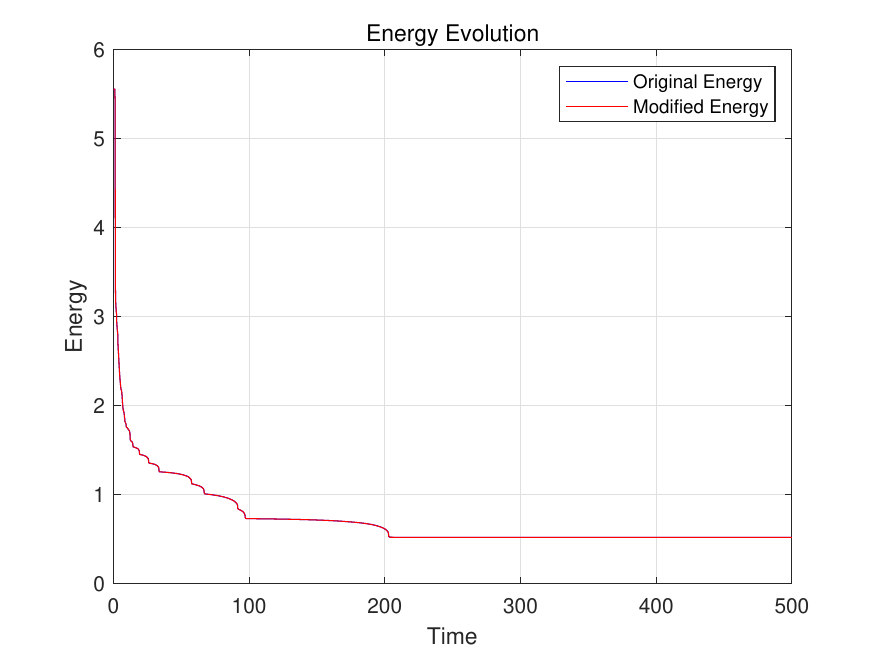}
}
\subfigure[R-IMEX-BDF2  with $\tau_n=1e-3$.]{
\includegraphics[width=2.1in]{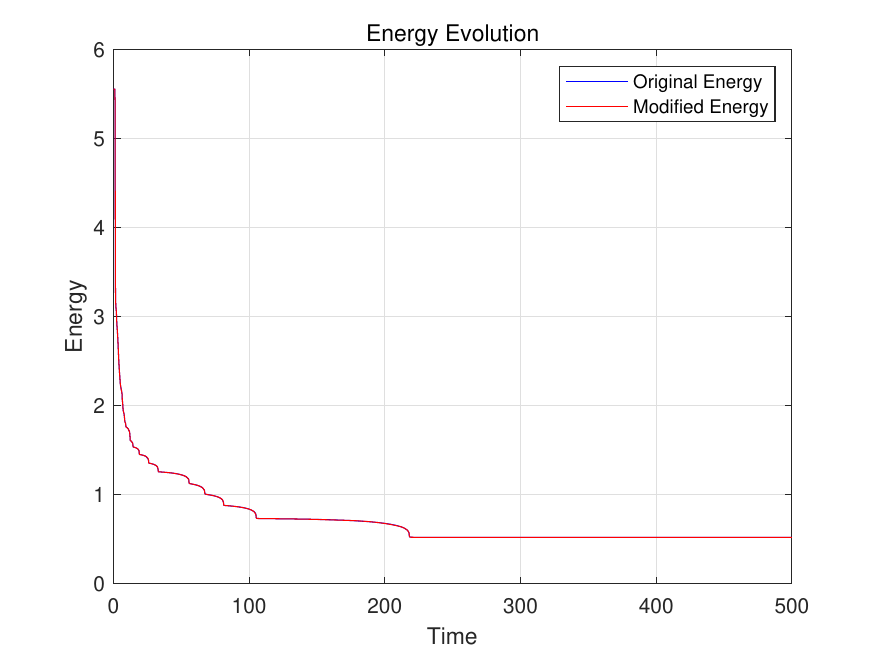}
}
\subfigure[Algorithm \ref{algo2} with (\ref{par6})(\ref{par10}).]{
\includegraphics[width=2.1in]{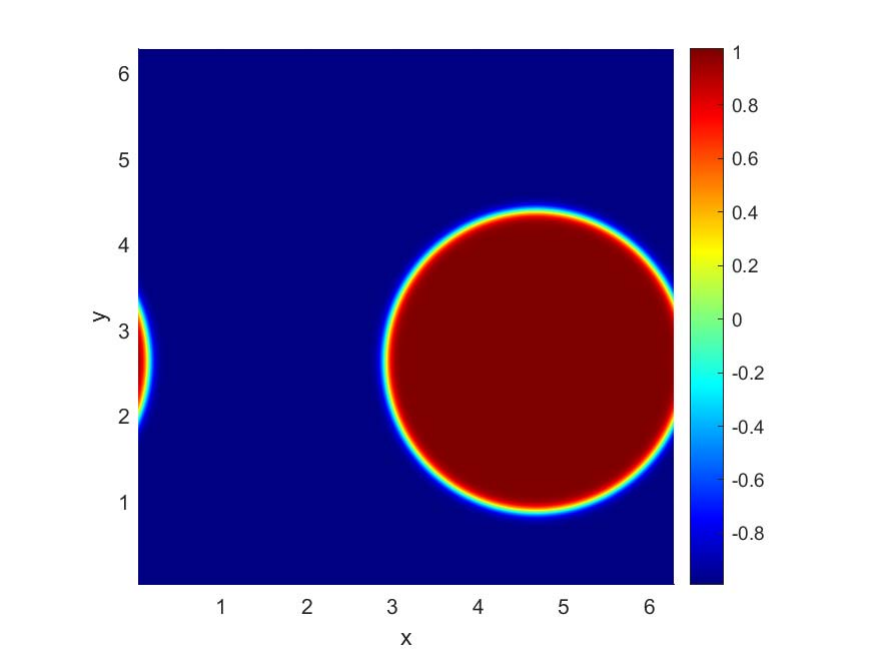}
}
\hspace{-8mm}
\subfigure[R-IMEX-BDF2A with (\ref{par3}).]{
\includegraphics[width=2.1in]{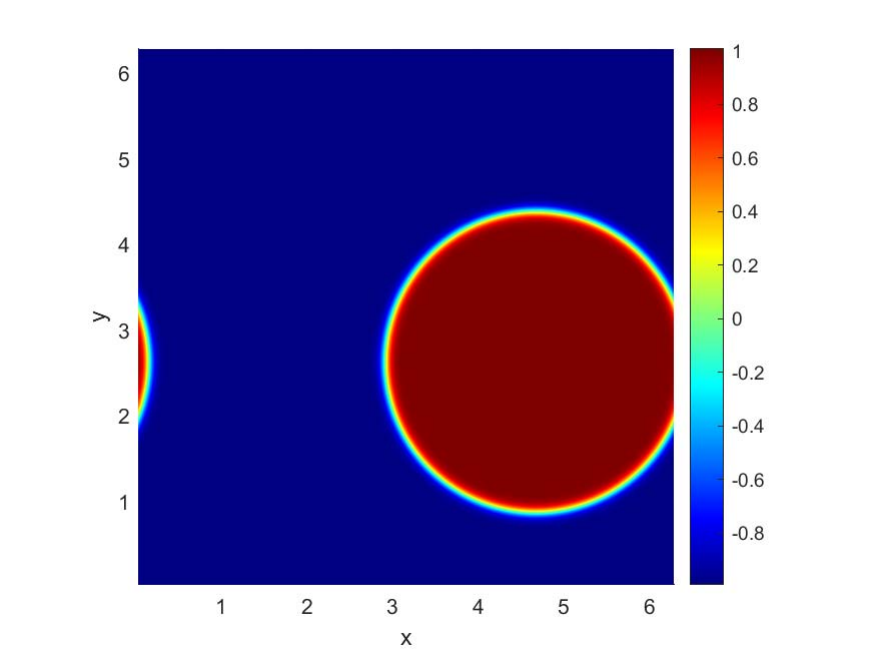}
}
\subfigure[R-IMEX-BDF2  with $\tau_n=1e-3$.]{
\includegraphics[width=2.1in]{T500_256.pdf}
}

\caption{(Example \ref{exp3})Modified energy and original energy history curves (a)(b)(c) of schemes with different parameters; The snapshots (d)(e)(f)  of the interfaces between two phases are depicted at $T=500$.}
\label{fig3}
\end{figure}

    \begin{table}[!ht]
		\centering
        \caption{(Example \ref{exp3})The  CPU time for some  algorithms at $T=500$.}
		\begin{tabular}{p{6cm} p{3cm} p{2cm} p{2cm}}
			\bottomrule
			Scheme/Algorithm &  CPU time &  Time \\
			\midrule
  Algorithm \ref{algo2} with  (\ref{par6})(\ref{par10})  & 1331.2s & \ \ \ 500 \\
    R-IMEX-BDF2A with (\ref{par3})  & 4342.9s & \ \ \ 500 \\
    R-IMEX-BDF2  with $\tau_n=1e-3$  & 12498.2s & \ \ \ 500 \\
  \bottomrule
\end{tabular}
		\label{tab7}
\end{table}
\par

Figures \ref{fig4} and \ref{fig5} show the histories of time steps and time step ratios for Algorithm \ref{algo2} (with parameters (\ref{par6}) and (\ref{par10})) and R-IMEX-BDF2A (with parameter (\ref{par3})), respectively, over the interval \((0,3)\). In Figure \ref{fig5}(a), R-IMEX-BDF2A with (\ref{par3}) has an average time step of \( 1.5 \times 10^{-3} \), while In Figure \ref{fig4}(a), the average time step for Algorithm \ref{algo2} is \( 2.9 \times 10^{-3} \). In Figures \ref{fig4}(b) and \ref{fig5}(b), the maximum time step ratios reach 8.8 and 5.2, respectively. As discussed in \cite{hou2023implicit, liao2022adaptive}, to maintain unconditional stability for variable time-stepping schemes, the step size ratio should not exceed 4.8645. Our approach, similar to the adaptive algorithm discussed in \cite{huang2020highly}, achieves unconditional energy stability without restricting the time-step ratio.

 \begin{figure}[!ht]
	\centering
	\subfigure[History of $\tau_n$ at (0,3).]{
		\includegraphics[width=7.3cm]{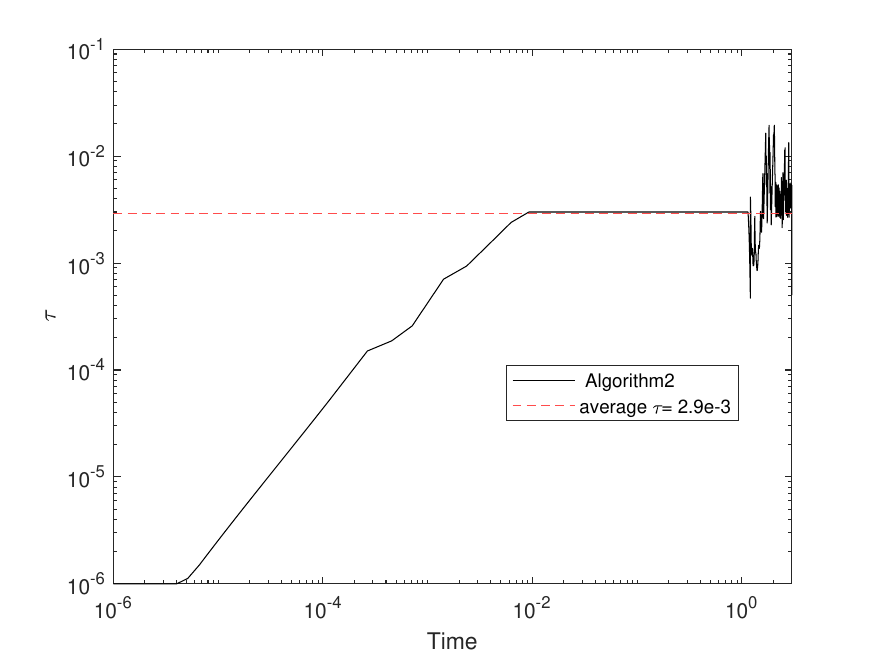}
	}
	\hspace{-10mm}
	\subfigure[History of $\tau_{n+1}/\tau_{n}$.]{
		\includegraphics[width=7.3cm]{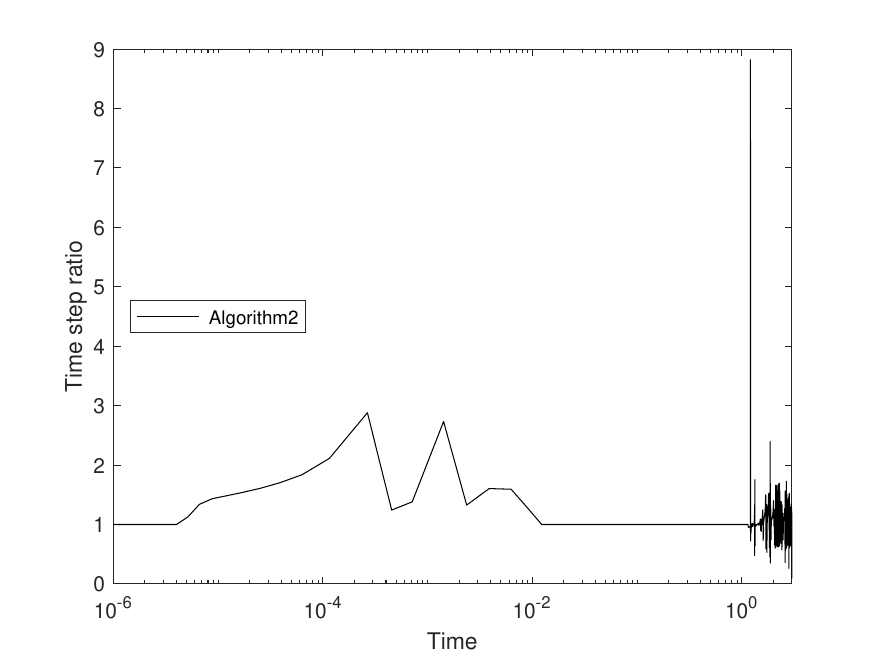}}
	\quad
	\caption{(Example \ref{exp3})Histories of time steps and time step ratios obtained by Algorithm \ref{algo2} with (\ref{par6})(\ref{par10}).The average time step is 2.9e-3.}\label{fig4}
\end{figure}
 \begin{figure}[!ht]
	\centering
	\subfigure[History of $\tau_n$  at (0,3).]{
	\includegraphics[width=7.3cm]{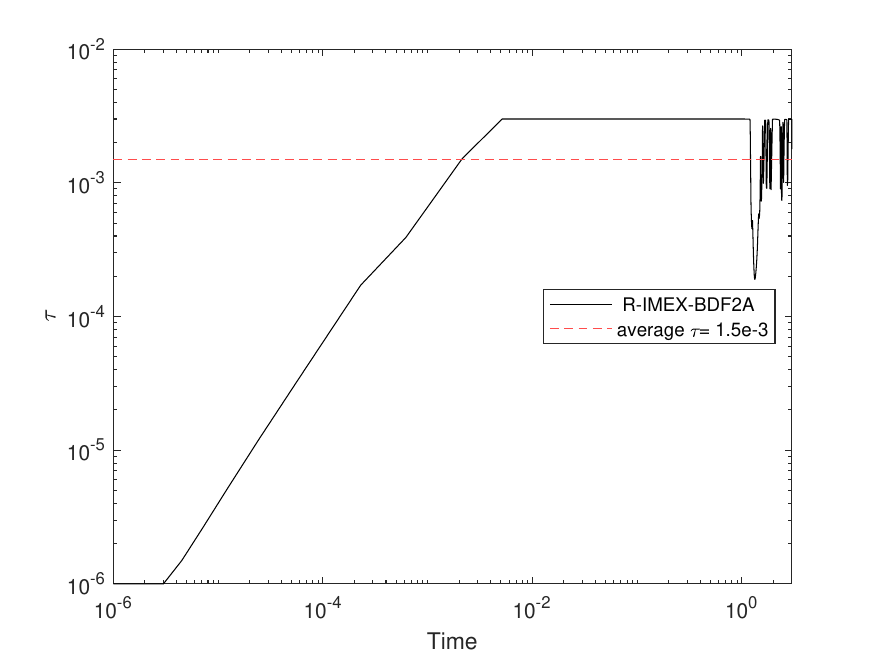}
}
\hspace{-10mm}
\subfigure[History of $\tau_{n+1}/\tau_{n}$.]{
	\includegraphics[width=7.3cm]{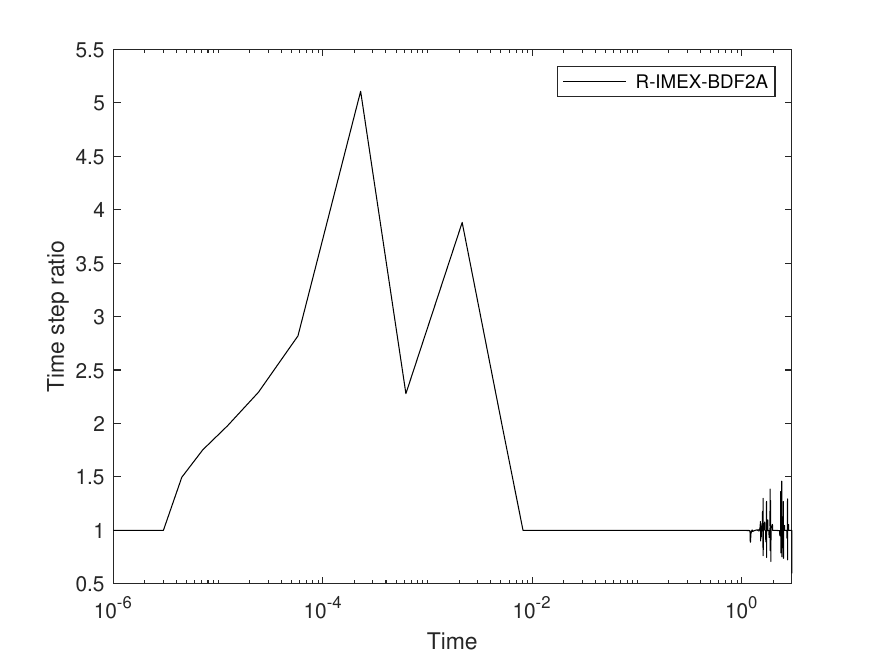}
}
\caption{(Example \ref{exp3})Histories of time steps and time step ratios obtained by R-IMEX-BDF2A with (\ref{par3}).The average time step is 1.5e-3.}
	\quad
\label{fig5}
\end{figure}
\end{example}	
\begin{example}(Buoyancy-driven flow) \label{exp4}
We consider in this last example the buoyancy-driven flow that two heavy fluid layers sandwich a light fluid layer \cite{zheng2022new,han2018second,jiang2024highly}.
This flow arises from fluid motion due to density differences caused by temperature or concentration variations. It occurs when buoyant forces dominate the flow behavior.
We adopt the flow equation as follows:
\begin{eqnarray*}
-\nabla\cdot[\nu(\phi)D(\u)]+\eta(\phi)\u=-\nabla p-\gamma\phi\nabla\mu+\textbf{b},\label{problem-1_4}
\end{eqnarray*}
where $\textbf{b}$ is a buoyancy term. We take the same assumption as in \cite{han2018second} that applying the following Boussinesq type approximation :
\begin{equation*}
\textbf{b}=(0,- {b(\phi)})^T, \ \ {b}(\phi)=\lambda(\phi-\bar{\phi}),
\end{equation*}
where $\bar{\phi}$ is the average of $\phi$, and $\lambda$ is the dimensionless parameter related to buoyancy strength \cite{le2006interfacial}.
In the simulation, we take the following  physical parameters:
\begin{equation*}
\epsilon=5e-2, \gamma=6, \eta=1, \nu=0.2,  \lambda=1.2,
\end{equation*}
and let the mobility that combines the Peclet number defined as:
\begin{equation*}\label{M}
	M(\phi)=\dfrac{1}{Pe}\sqrt{(1+\phi)^{2}(1-\phi)^{2}+\epsilon^{2}},
\end{equation*}
where $Pe=1$. The initial values are set as
\begin{equation*}\label{tuophi0}
	\phi(x,y,0)=tanh(\dfrac{y-(0.5\pi-(0.5+0.1cos(x)))}{\sqrt{2}\epsilon})tanh(\dfrac{y-(0.5\pi+(0.5+0.1cos(x)))}{\sqrt{2}\epsilon}).
\end{equation*}
In this example, we set $T_c=0.4$ for Algorithm \ref{algo2} and take the following numerical  parameters as follows:
\begin{eqnarray}\label{par8}
&\text{R-IMEX-BDF3A}:\rho=0.75,tol=1e-3,r=0.33,m=1,\gamma^{*}=1200,\tau_{max}=1e-3;	\\
 \label{par9}
&\text{R-IMEX-BDF2A}:\rho=0.89,tol=1e-2,r=0.30,m=1,\gamma^{*}=1200,\tau_{max}=1e-2.	
\end{eqnarray}

 Snapshots of the phase field evolution at times \( T = 0, 3, 5, 7, 7.8 \) are shown in Figure \ref{fig7}. These snapshots clearly illustrate the buoyancy-driven flow physics we simulate. The blue fluid becomes thicker on the sides and thinner in the middle, as it is lighter than the red fluid. Due to buoyancy strength, the blue fluid on both sides rises, while the fluid in the middle moves downward. Over time, the blue fluid in the middle breaks apart and forms floating spheroids. This result is in good agreement with simulations in \cite{zheng2022new,han2018second,jiang2024highly}, demonstrating the capability of the proposed scheme to capture complex phenomena.
	\begin{figure}[!ht]
	\centering
	\subfigure[$T=0$]{
		\begin{minipage}[t]{0.18\linewidth}
			\includegraphics[width=1.5in]{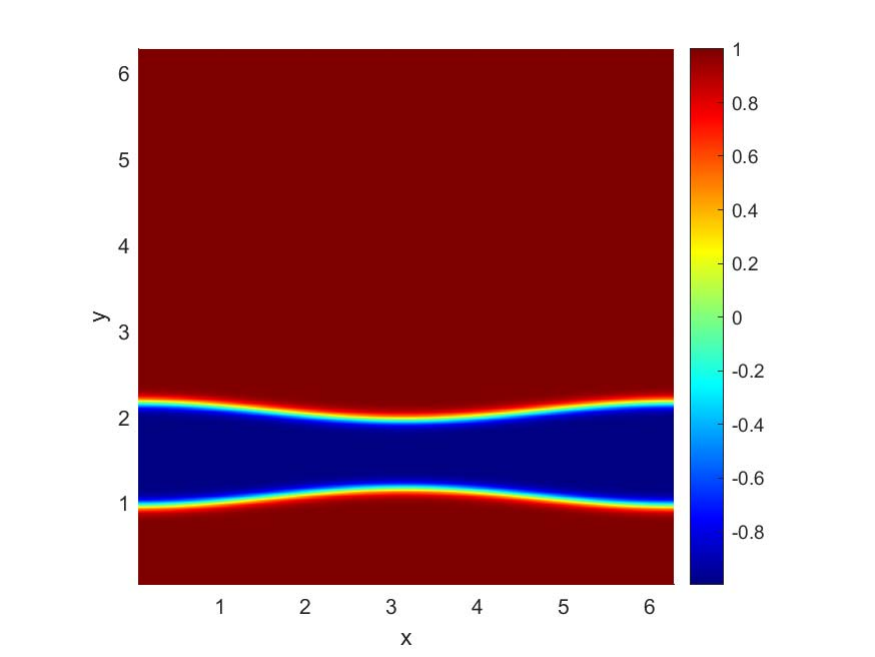}
		\end{minipage}%
	}%
	\subfigure[$T=3$]{
		\begin{minipage}[t]{0.18\linewidth}
			\includegraphics[width=1.5in]{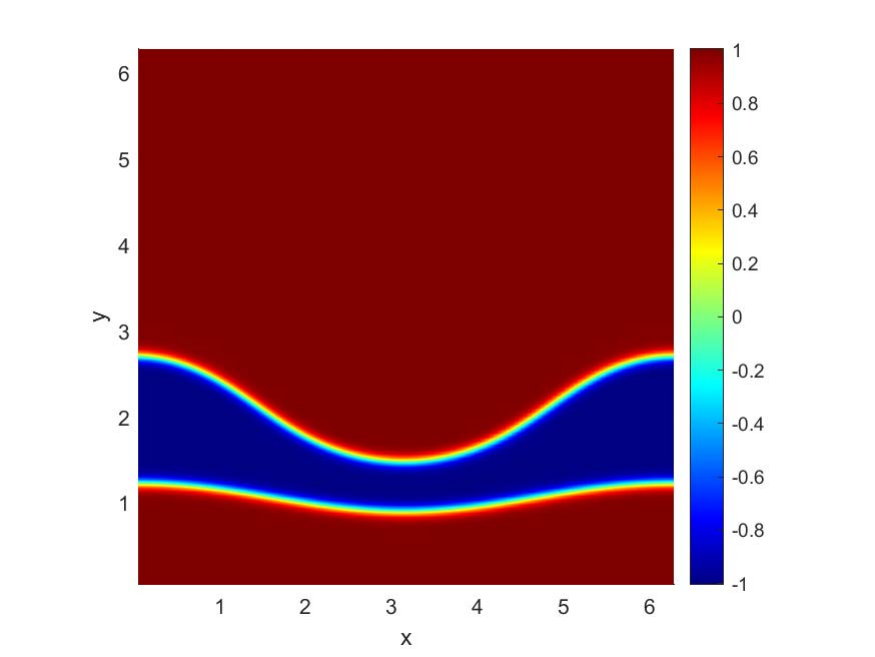}
		\end{minipage}%
	}%
	\subfigure[$T=5$]{
		\begin{minipage}[t]{0.18\linewidth}
			\includegraphics[width=1.5in]{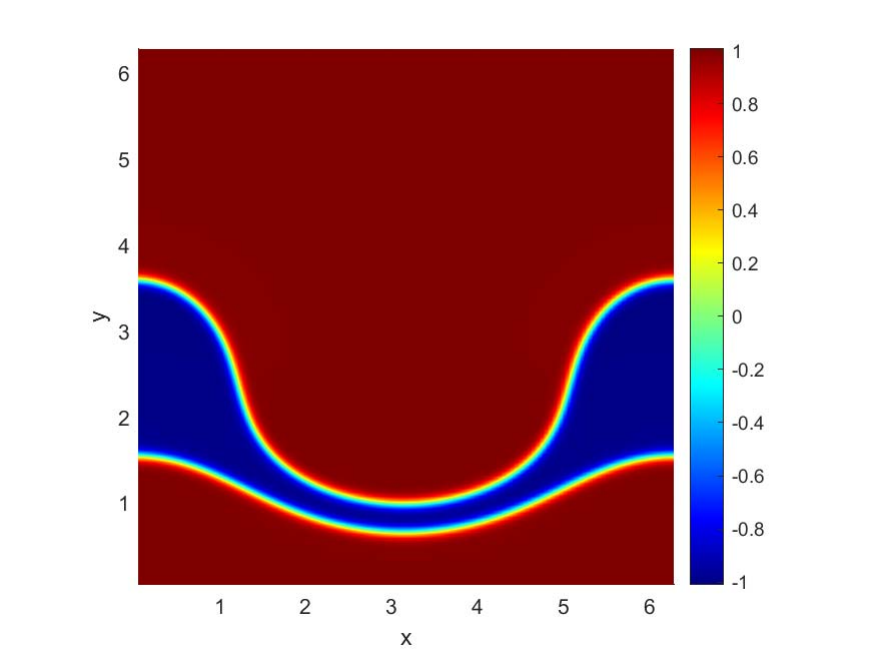}
		\end{minipage}%
	}%
	\subfigure[$T=7$]{
		\begin{minipage}[t]{0.18\linewidth}
			\includegraphics[width=1.5in]{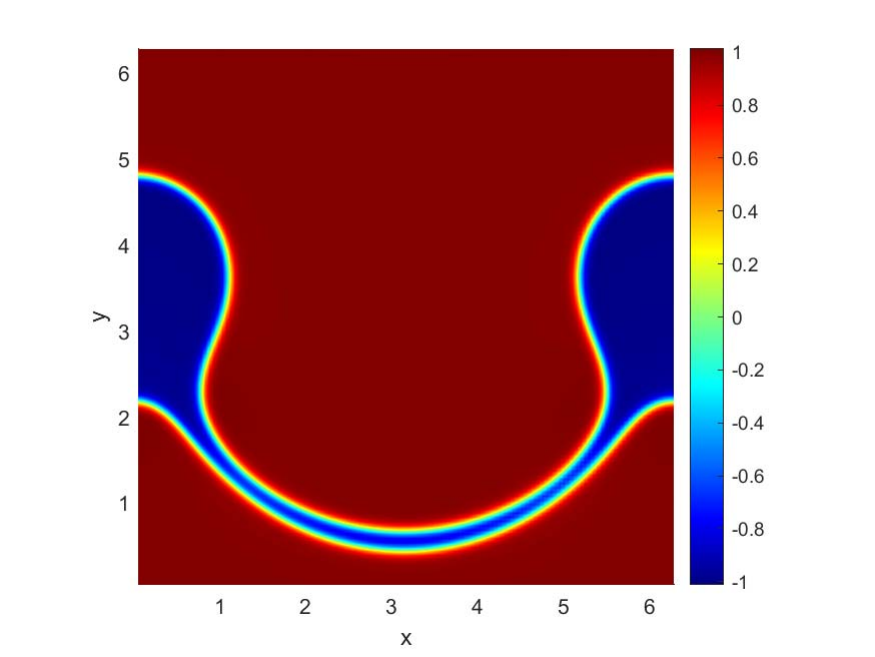}
		\end{minipage}%
	}%
	\subfigure[$T=7.8$]{
		\begin{minipage}[t]{0.18\linewidth}
			\includegraphics[width=1.5in]{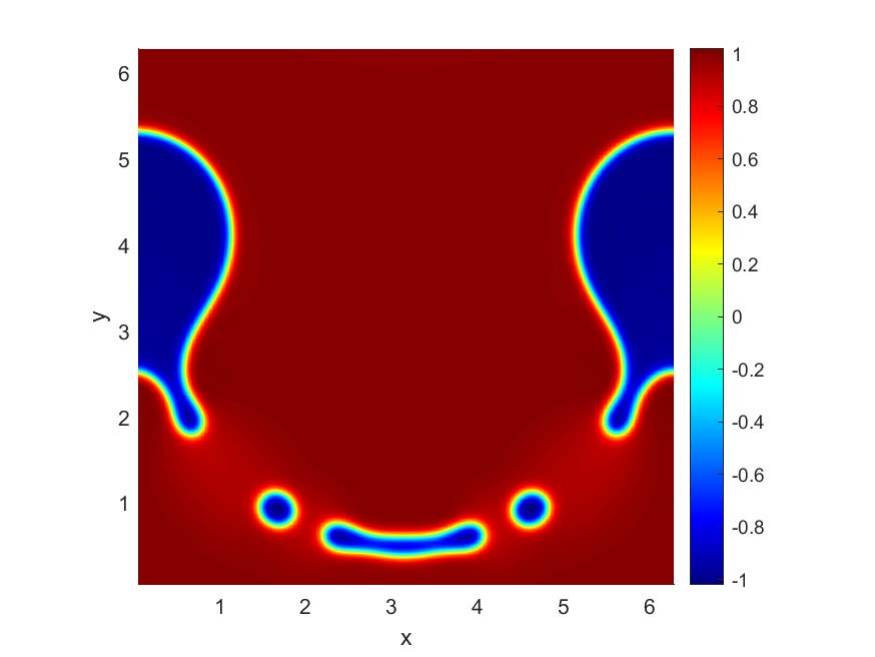}
		\end{minipage}%
	}%
	\centering
	\caption{(Example \ref{exp4})Snapshots of  buoyancy-driven flow  function $\phi$ by R-IMEX-BDF2A with (\ref{par9}).}
	\label{fig7}
\end{figure}	
\par
Table \ref{tab8} further corroborates the findings in Example \ref{exp3}, highlighting that Algorithm \ref{algo2} leads to a substantial improvement in accuracy compared to the R-IMEX-BDF2A scheme. The adaptive time-stepping strategy proposed here reduces computational cost or enhances accuracy compared to the  R-IMEX-BDF2 or IMEX-BDF2 schemes with fixed time steps. Additionally, with a fixed time step of $\tau_n=2e-3$, the R-IMEX-BDF2 scheme improves accuracy compared to the IMEX-BDF2 scheme, consistent with the results in Example \ref{exp2}. Figure \ref{fig8} presents similar numerical results from different algorithms at \( T=7.8 \).
It is important to note that numerous experiments show that strictly following the update time step strategy from \cite{huang2020highly} or \cite{zhang2013adaptive, hou2023implicit, cheng2018multiple}, such as setting $\gamma^{\star}_k = 0$ or using a large $tol_k$, will lead to higher computational costs or lower accuracy compared to the proposed methods herein.
This indicates that proper error indicators and energy-based step size control help in designing an adaptive time-stepping strategy. They enable adaptive time-stepping methods that simulate complex phenomena more efficiently.

\par
\begin{table}[!ht]
	\centering
    \caption{(Example \ref{exp4})The $L^{2}$-norm error of $\phi$ and CPU time for some  algorithms at $T=7.8$.}
	\begin{tabular}{p{6cm} p{3cm} p{2cm} p{2cm}}
		\bottomrule
		Scheme/Algorithm & $L^{2}$-norm error & CPU time & Time $T$\\
		\midrule
		Algorithm \ref{algo2} with (\ref{par8})(\ref{par9}) & 9.331e-3 & 578.6s & \ \ \ 7.8 \\
		R-IMEX-BDF2A with (\ref{par9}) & 6.142e-2 & 700.8s & \ \ \ 7.8\\
        R-IMEX-BDF2 with $\tau_n=2e-3$ & 1.053e-2 & 2032.4s &  \ \ \ 7.8\\
        IMEX-BDF2 with $\tau_n=2e-3$ & 1.743e-2 & 1950.0s &  \ \ \ 7.8\\
		\bottomrule
	\end{tabular}
	\label{tab8}
\end{table}
\end{example}
\begin{figure}[!ht]
\centering
\subfigure[Algorithm \ref{algo2} with (\ref{par8})(\ref{par9}).]{
\includegraphics[width=2.1in]{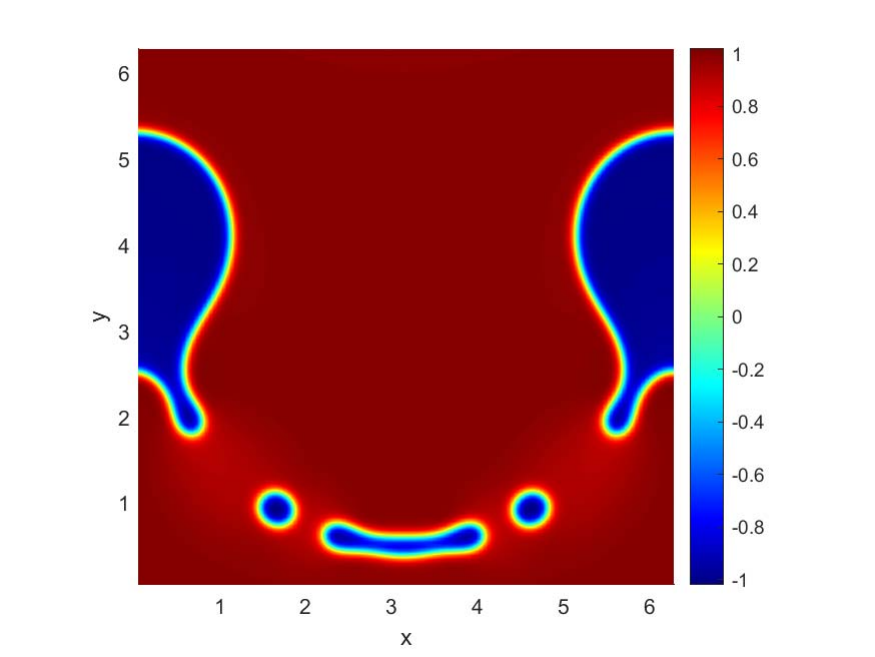}
}
\hspace{-8mm}
\subfigure[R-IMEX-BDF2 with $\tau_n=0.002$.]{
\includegraphics[width=2.1in]{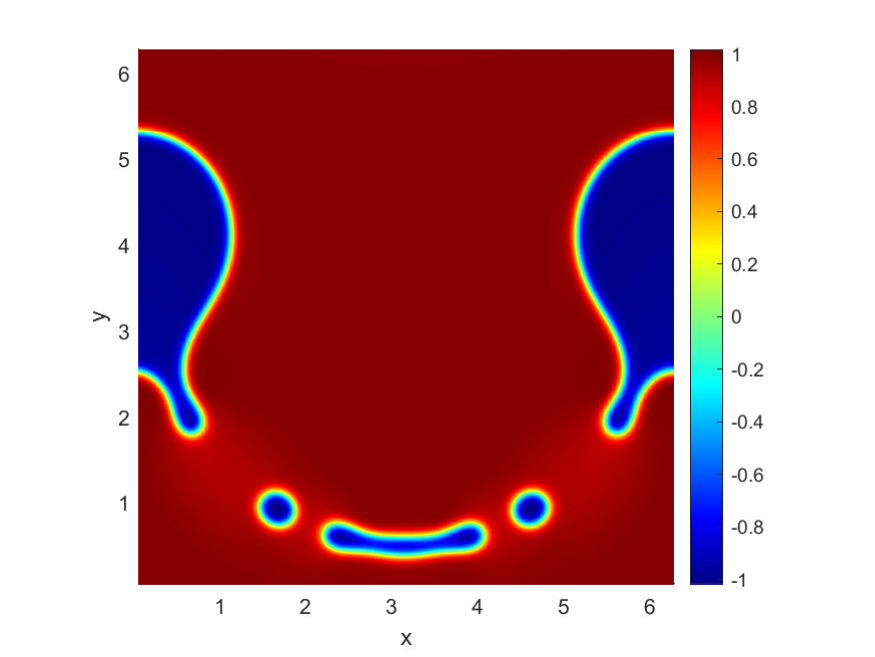}
}
\subfigure[IMEX-BDF2 with $\tau_n=0.002$.]{
\includegraphics[width=2.1in]{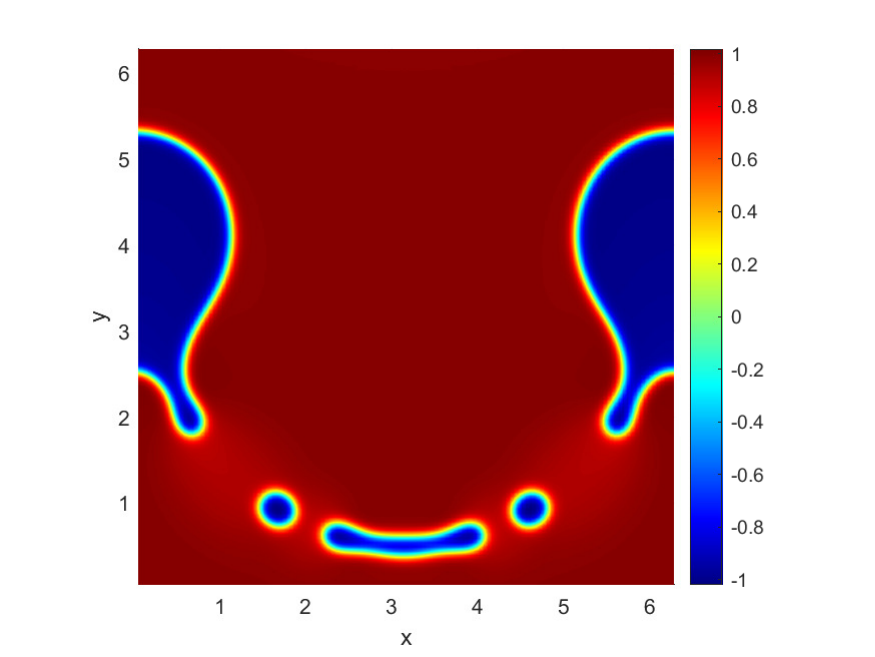}
}
\caption{(Example \ref{exp4})The snapshots of  buoyancy-driven flow  function $\phi$ by different algorithms at T=7.8.}\label{fig8}
\end{figure}

\section{Concluding  remarks}
For the Cahn-Hilliard-Brinkman system, we proposed unconditionally energy-stable R-IMEX-BDF$k$ schemes. These schemes only require solving elliptic and algebraic equations at each step. Compared to IMEX-BDF$k$ schemes, the relaxation improves accuracy.
Furthermore, We utilized the properties of the R-IMEX-BDF$k$ scheme and the coupled system to design effective error indicators and an energy-based time-step control strategy. This led to the development of a more flexible adaptive strategy, the R-IMEX-BDF$k$A algorithm. Then,
based on the properties of high-order  BDF schemes and the CHB model, we developed a hybrid-order adaptive algorithm. In comparison with fixed time step methods, this adaptive algorithm significantly improves computational efficiency.
It employs R-IMEX-BDF$3$A to capture the rapid phase evolution in the early stage and R-IMEX-BDF$2$A for long-term dynamics, enhancing overall computational efficiency. Numerical experiments validated the accuracy and effectiveness of the proposed methods. Notably, these algorithms can be directly extended to other phase-field equations. Future work will focus on developing similar adaptive time-stepping strategies for different systems, providing a promising avenue for further research.
\bibliographystyle{siam}
\bibliography{ref}
\end{document}